\documentclass[11pt,a4paper]{article}

\usepackage{amsmath}
\usepackage{amsfonts}
\usepackage{amssymb}
\usepackage{algorithm}
\usepackage{algorithmic}
\usepackage{color}
\usepackage{enumerate}
\usepackage{graphicx}
\usepackage[hidelinks]{hyperref}
\usepackage{lscape}
\usepackage{longtable}
\usepackage{mathrsfs}
\usepackage{rotating}
\usepackage{subcaption}
\usepackage{url}
\usepackage{titlesec} 
\usepackage[titletoc]{appendix}

\newtheorem{theorem}{Theorem}[section]

\newtheorem{corollary}{Corollary}[section]

\newtheorem{lemma}{Lemma}[section]

\newtheorem{assumption}{Assumption}[section]

\newenvironment{proof}[1][Proof]{\textbf{#1.} }{\hfill$\Box$}

\def\flow{f_{\mathrm{low}}}

\def\R{\mathbb{R}}
\def\N{\mathbb{N}}
\def\1st{\mathrm{G}}
\def\2nd{\mathrm{H}}

\def\cO{\mathcal{O}}
\def\T{\mathrm{T}}
\def\Ir{\mathcal{R}}
\def\It{\mathcal{N}}

\newcommand{\revised}[1]{\textcolor{black}{#1}}

\begin{document}

\title{A nonlinear conjugate gradient method with complexity guarantees 
%\revised{for nonconvex optimization}}
and its application to nonconvex regression} 
\author{R\'emi Chan{-}{-}Renous-Legoubin
\thanks{Universit\'e Paris Dauphine-PSL, 75016 Paris, France 
({\tt remi.chan-renous-legoubin@dauphine.eu}).}
\and Cl\'ement W. Royer
\thanks{LAMSADE, CNRS, Universit\'e Paris Dauphine-PSL, 75016 Paris, France
({\tt clement.royer@lamsade.dauphine.fr}). 
Support for this research was provided by CNRS INS2I under the grant GASCON 
and by Agence Nationale de la Recherche through program ANR-19-P3IA-0001 (PRAIRIE 3IA Institute)
.}
}
%\date{}
\maketitle
\thispagestyle{empty}
{\small
\begin{abstract}
	Nonlinear conjugate gradients are among the most popular
	techniques for solving continuous optimization problems. Although 
	these schemes have long been studied from a global convergence 
	standpoint, their worst-case complexity properties have yet to be 
	fully understood, especially in the nonconvex setting. In particular, 
	it is unclear whether \revised{nonlinear conjugate gradient} methods possess better 
	guarantees than first-order methods such as gradient descent. \revised{Meanwhile,
	recent experiments} have shown \revised{impressive} performance of 
	standard nonlinear conjugate gradient \revised{techniques} on \revised{certain} 
	nonconvex problems, 
	even when compared with methods endowed with the best known complexity 
	guarantees. 	
	
	In this paper, we propose a nonlinear conjugate gradient \revised{scheme} 
	based on a simple line-search paradigm and a modified restart condition.
	These two ingredients allow for monitoring the 
	properties of the search directions, which is instrumental 
	in obtaining complexity guarantees. Our complexity results illustrate 
	the possible discrepancy between nonlinear conjugate gradient methods 
	and classical gradient descent. \revised{A numerical investigation on 
	nonconvex robust regression problems as well as a standard benchmark 
	illustrate that the restarting 
	condition can track the behavior of a standard implementation.}
\end{abstract} 
}

%\tableofcontents

%\newpage

%%%%%%%%%%%%%%%%%%%%%%%%%%%%%%%%%%%%%%%%%%%%%%%%%%%%%%%%%%%%%%%%%%%%%%%%%%%%%%%
\section{Introduction}
\label{sec:intro}

In this paper, we are interested in the problem
\begin{equation}
\label{eq:pb}
	\min_{x \in \R^n} f(x),
\end{equation}
where $f: \R^n \rightarrow \R$ is a Lipschitz continuously differentiable, 
nonconvex function.
Given a tolerance $\epsilon \in (0,1)$, our goal is to compute 
an $\epsilon$-approximate stationary point, that is, a vector $x \in \R^n$ such that
\begin{equation}
\label{eq:epspoint}
	\|\nabla f(x)\| \le \epsilon,
\end{equation}
where $\|\cdot\|$ denotes the Euclidean norm in $\R^n$.
We are interested in algorithms that provably reach a point 
satisfying~\eqref{eq:epspoint} in a finite amount of computational 
work (number of iterations, and evaluations of $f$ and its derivatives). 
Such results are termed \emph{worst-case complexity guarantees}, and have been a 
growing area of 
interest in nonconvex optimization, fueled by the interest of such results in 
data analysis and statistics~\cite{SJWright_2018}. 

The most classical algorithm for solving problem~\eqref{eq:pb} is gradient 
descent, that proceeds by moving along the negative gradient direction. It is 
known that such a method typically reaches a point satisfying~\eqref{eq:epspoint} 
in at most $\cO(\epsilon^{-2})$ iterations or gradient evaluations\footnote{Throughout this paper, we write $\cO(A)$ to indicate a constant times $A$, where the terms in the constants do not 
depend on that in $A$. The notation $\tilde{\cO}(A)$ stands for 
$\cO(\log^c(A) A)$ with $c>0$.}. This bound contrasts with its counterparts in convex 
and strongly convex problems, and is sharp for gradient descent, even when exact line searches are used~\cite{CCartis_NIMGould_PhLToint_2010,
CCartis_PhRSampaio_PhLToint_2015}. 
\revised{Perhaps more suprisingly, the bound $\cO(\epsilon^{-2})$ is 
also sharp for more elaborate frameworks, such as trust region and 
Newton's method~\cite{CCartis_NIMGould_PhLToint_2010}.}
Under additional regularity assumptions 
on the objective function, it is possible to design algorithms 
with better properties, that typically require higher-order information 
to be used explicitly in the algorithm. \revised{Cubic regularization methods, that achieve an iteration complexity 
in $\cO(\epsilon^{-3/2})$, were the first class of algorithms to be 
equipped with such a guarantee in the nonconvex setting~\cite{YuNesterov_BTPolyak_2006,CCartis_NIMGould_PhLToint_2011b}. 
These results attracted quite a lot of attention, leading to multiple second-order 
methods with similar complexity bounds being proposed~\cite{EGBirgin_JMMartinez_2017, CCartis_NIMGould_PhLToint_2011c,CCartis_NIMGould_PhLToint_2019a,FECurtis_DPRobinson_MSamadi_2017,FECurtis_DPRobinson_MSamadi_2019,
CWRoyer_SJWright_2018,PhLToint_2013}. Although out of the scope of this work, we note 
that the complexity results can be further improved by leveraging high-order 
smoothness~\cite{EGBirgin_JLGardenghi_JMMartinez_SASantos_PhLToint_2017}.}

\revised{More recently,} accelerated gradient techniques combined with 
negative curvature estimation procedures, for which a gradient evaluation 
complexity in $\tilde{\cO}(\epsilon^{-7/4})$ can be shown~\cite{YCarmon_JCDuchi_OHinder_ASidford_2017a,
YCarmon_JCDuchi_OHinder_ASidford_2018}. \revised{The latter methods depart from 
standard nonlinear optimization techniques, and were shown to be outperformed 
by a standard nonlinear conjugate gradient implementation on a nonconvex 
regression example~\cite{YCarmon_JCDuchi_OHinder_ASidford_2017a}. However, 
their construction facilitates the derivation of complexity results.}
Providing a complexity analysis of standard nonlinear optimization 
schemes, on the other hand, remains a challenging endeavor.

Conjugate gradient (CG) methods are an example of popular nonlinear 
optimization techniques that have yet to be endowed with a comprehensive 
complexity analysis. When applied to strongly convex quadratics, those 
methods typically reduce to linear conjugate gradient, and complexity 
guarantees have long been known in that setting~\cite{JNocedal_SJWright_2006}. 
Recent results have shown that 
linear conjugate gradient can also be used in conjunction with a 
Newton-type framework when the objective is Lipschitz twice 
continuously differentiable. The resulting methods 
can be analyzed from a complexity viewpoint on nonconvex optimization 
problems, leading to complexity bounds that match the best known 
in the literature for this class of problems~\cite{FECurtis_DPRobinson_CWRoyer_SJWright_2021, 
CWRoyer_MONeill_SJWright_2020,CWRoyer_SJWright_2018}.

The situation becomes quite different when considering nonlinear conjugate 
gradient techniques. Numerous conditions for global convergence of nonlinear CG 
schemes have been proposed in the literature~\cite{YHDai_2002,
LGrippo_SLucidi_1997,LGrippo_SLucidi_2005,WWHager_HZhang_2005,
WWHager_HZhang_2006a,WWHager_HZhang_2006b}, yet early complexity 
analyzes showed that nonlinear CG techniques could have worse guarantees than gradient descent on strongly convex problems~\cite{ASNemirovski_DBYudin_1983}. 
Still, recent proposals have combined modern nonlinear conjugate gradient with accelerated gradient 
tools to yield a method that both reduces to linear CG on quadratics 
and converges as fast as accelerated gradient on convex and strongly 
convex problems~\cite{SKarimi_SAVavasis_2016,SKarimi_SAVavasis_2017,
SKarimi_SAVavasis_2021}. Nevertheless, deriving complexity guarantees 
for classical nonlinear CG variants remains a difficult task, particularly 
in the nonconvex setting. 

In this paper, we propose a nonlinear conjugate gradient \revised{framework for 
optimizing nonconvex functions with complexity guarantees}. By adapting the classical restart condition 
of nonlinear conjugate gradient methods, we are able to provide decrease guarantees 
at every iteration, that depend on whether the restart condition is triggered. 
Overall, our method is shown to depart from gradient descent on certain 
iterations, for which our analysis provides better guarantees. \revised{To 
investigate the practical impact of this condition, we first confirm 
earlier findings about the performance of these 
methods on nonconvex regression tasks, where regimes in which the restart condition is representative 
of the algorithmic behavior of nonlinear conjugate gradient can be identified. 
We also investigate the difference between classical implementations  of 
nonlinear CG that lack complexity guarantees and our framework on a nonlinear optimization 
benchmark.}

The rest of this paper is organized as follows. In Section~\ref{sec:ncg}, we 
recall the key features of nonlinear conjugate gradient algorithms, and we 
describe our framework based on a modified restart condition. Complexity 
results for this framework are obtained and discussed in 
Section~\ref{sec:wcc}. \revised{In Section~\ref{sec:num}, we then conduct 
a numerical study of our proposed method involving both nonconvex regression 
tasks and comparisons on a nonlinear optimization test set.}
Final comments are made in Section~\ref{sec:conc}.

%\paragraph{Notations} Throughout this paper, we use $\cO(\cdot)$ to denote

%%%%%%%%%%%%%%%%%%%%%%%%%%%%%%%%%%%%%%%%%%%%%%%%%%%%%%%%%%%%%%%%%%%%%%%%%%%%%%%
\section{Nonlinear conjugate gradient framework}
\label{sec:ncg}
%%%%%%%%%%%%%%%%%%%%%%%%%%%%%%%%%%%%%%%%%%%%%%%%%%%%%%%%%%%%%%%%%%%%%%%%%%%%%%%

In this section, we describe a nonlinear conjugate gradient method based on 
Armijo line-search and a modified restart condition. To this end, we first 
recall the main features of nonlinear conjugate gradient methods, then 
provide a description of our proposed scheme.

%%%%%%%%%%%%%%%%%%%%%%%%%%%%%%%%%%%%%%%%%%%%%%%%%%%%%%%%%%%%%%%%%%%%%%%%%%%%%%%
\subsection{Nonlinear conjugate gradient}
\label{subsec:pb:ncgclassic}

Nonlinear conjugate gradient techniques~\cite{WWHager_HZhang_2006b,
RPytlak_2009} are iterative optimization schemes of the form
\begin{equation}
\label{eq:ncgit}
	x_{k+1}=x_k+\alpha_k d_k,
\end{equation}
where the direction $d_k$ is a search direction, and $\alpha_k$ is a 
stepsize typically computed through a line search. Several line-search 
strategies have been proposed in the literature~\cite{YHDai_2002,WWHager_HZhang_2006b}. 
In this paper, we focus on using a backtracking 
Armijo line search, \revised{since this line search led to
good performance on nonconvex regression problems}~\cite{YCarmon_JCDuchi_OHinder_ASidford_2017a}.

The goal 
of a conjugate gradient approach is to combine local information (i.e. the 
negative gradient at the current point) with the \emph{previous direction}, 
that is possibly still relevant for the next iterate. As a result, a 
nonlinear CG method selects $d_0 = -\nabla f(x_0)$ and
\begin{equation}
\label{eq:ncgdk}
	\left\{
		\begin{array}{ll}
			d_0 = -\nabla f(x_0) & \\
    		d_{k+1} = -\nabla f(x_{k+1}) + \beta_{k+1} d_k 
    		&\forall k \in \mathbb{N}.
    	\end{array}
    \right.
\end{equation}
The choice of the formula for \revised{the parameter} $\beta_{k+1}$ gives 
rise to various nonlinear CG schemes~\cite{WWHager_HZhang_2006a}. The most 
popular choices for $\beta_{k+1}$ include the \emph{Fletcher-Reeves} formula
\begin{equation}
\label{eq:ncgFR}
    \beta_{k+1}^{FR} := \frac{\|\nabla f(x_{k+1})\|^2}{\|\nabla f(x_{k})\|^2}
\end{equation}
and the \emph{Polak-Ribi\`ere} (also known as Polak-Ribi\`ere-Polyak) formula
\begin{equation}
\label{eq:ncgPR}
    \beta_{k+1}^{PR} := \frac{\nabla f(x_{k+1})^T (\nabla f(x_{k+1})-\nabla f(x_{k}))}{\|\nabla f(x_{k})\|^2}.
\end{equation}

%A similar reasoning can be adopted when the formula for $\beta_{k+1}$ yields 
%a nonpositive value. 
A popular variant of the Polak-Ribi\`ere formula is the
\emph{PRP+} update:
\begin{equation}
\label{eq:ncgPRP+}
    \beta_{k+1}^{PRP+} := \max\left\{
    \frac{\nabla f(x_{k+1})^T (\nabla f(x_{k+1})-\nabla f(x_{k}))}
    {\|\nabla f(x_{k})\|^2},0\right\}.
\end{equation}
Using the PRP+ formula has been shown to guarantee global convergence of a 
nonlinear conjugate gradient method with appropriate line-search 
condition~\cite{JCGilbert_JNocedal_1992}. 
\revised{More recently, the formula proposed by Hager and 
Zhang~\cite{WWHager_HZhang_2006a} 
\begin{equation}
\label{eq:ncgHZ}
	\beta_{k+1}^{HZ} := \frac{1}{d_{k}^\T y_{k}}\left( y_{k} 
	- 2 d_{k} \frac{\|y_{k}\|^2}{d_{k}^\T y_{k}} \right)^\T \nabla f(x_{k+1}),
\end{equation}
where $y_k = \nabla f(x_{k+1}) - \nabla f(x_k)$, has been found quite 
successful in practice.
}

\revised{Regardless of the variant that is used, a nonlinear conjugate 
gradient method can produce iterates such that} $\nabla f(x_{k+1})^T d_{k+1} \ge 0$, in 
which case there is no guarantee for decrease in the direction $d_{k+1}$: 
a typical fix consists in redefining the search direction as 
$d_{k+1}=-\nabla f(x_{k+1})$. This process is called restarting, 
\revised{and is a common feature of a nonlinear conjugate gradient 
implementation. In the next section, we propose an alternative to this 
restart condition that endows a conjugate gradient scheme with 
complexity results.}
%Following previous work and 
%motivated by this restarting aspect of the PRP+ update, we will adopt 
%formula~\eqref{eq:ncgPRP+} in our algorithm. Our method also exploits the 
%restart condition, as described in the next section.

%%%%%%%%%%%%%%%%%%%%%%%%%%%%%%%%%%%%%%%%%%%%%%%%%%%%%%%%%%%%%%%%%%%%%%%%%%%%%%%
\subsection{Nonlinear CG with modified restart}
\label{subsec:pb:algo}

In our framework, we build on the restarting idea by monitoring the 
value of $\nabla f(x_k)^\T d_k$ and that of $\|d_k\|$. Algorithm~\ref{alg:algo} 
describes our framework. At every iteration, we perform a 
backtracking line search to compute a step that yields a suitable 
decrease in the objective function (see condition~\eqref{eq:lscond}).
Once the new point has been computed, we evaluate the gradient at 
the next iterate, as well as the parameter $\beta_{k+1}$\revised{, that is 
typically chosen from one of the formulas given in the previous section. Both 
the gradient and the parameter are then used to define the new search direction.}
%, which we compute through the PRP+ update~\eqref{eq:ncgPRP+}. 

\begin{algorithm}[ht!]
\caption{Nonlinear conjugate gradient with modified restart condition}
\label{alg:algo}
\begin{algorithmic}
	\STATE \emph{Inputs:} $x_0 \in \mathbb{R}^d$, $\eta \in (0,1)$, $\theta \in (0,1)$, 
	$\sigma \in (0,1]$, 
	$\kappa \ge 1$, $p \ge 0$, $q \ge 0$.
	\STATE Set $g_0=\nabla f(x_0)$, $d_0 = -g_0$, $k = 0$.
	\FOR{$k=0,1,2,\dots$}
  		\STATE Compute $\alpha_k= \theta^{j_k}$ where $j_k$ is the smallest 
  		nonnegative integer such that
  		\begin{equation}
  		\label{eq:lscond}
  			f(x_k+\alpha_k d_k) < f(x_k) + \eta \alpha_k g_k^T d_k.
  		\end{equation}
  		\STATE Set $x_{k+1} = x_k + \alpha_k d_k$ and $g_{k+1} = \nabla f(x_{k+1})$.
  		\STATE \revised{Choose a conjugate direction parameter $\beta_{k+1}$.}
%  		\STATE Set $\beta_{k+1} = \max\left\{\tfrac{g_{k+1}^\T (g_{k+1}-g_k)}{\|g_k\|^2},0\right\}$.
  		\STATE Set $d_{k+1} = -g_{k+1} + \beta_{k+1} d_k$.
  		\STATE If the condition
  		\begin{equation}
  		\label{eq:restartcond}
  			g_{k+1}^\T d_{k+1} \ge -\sigma \|g_{k+1}\|^{1+p} 
  			\quad \mbox{or} \quad 
  			\|d_{k+1}\| \ge \kappa\|g_{k+1}\|^q,
  		\end{equation}
  		holds, restart the algorithm by setting $d_{k+1}=-g_{k+1}$.
  	\ENDFOR
\end{algorithmic}
\end{algorithm}

The key ingredient to Algorithm~\ref{alg:algo} is the restarting 
condition~\eqref{eq:restartcond}\revised{, that determines whether the 
nonlinear CG direction is kept for the next iteration.}
For any $k \ge 0$, if iteration $k$ does not end with a restart, we 
have\footnote{Although the inequalities should be strict, we use non-strict 
inequalities for notational convenience.}
\begin{equation}
\label{eq:conddk}
	g_{k+1}^\T d_{k+1} \le -\sigma \|g_{k+1}\|^{1+p}
	\quad \mbox{and} \quad 
	\|d_{k+1}\| \le \kappa \|g_{k+1}\|^q.
\end{equation}
\revised{In Section~\ref{sec:wcc}, we will establish that these inequalities 
together with the line search allow for proving complexity 
bounds for Algorithm~\ref{alg:algo}. Note that the above properties can be viewed as a 
more general case of the following condition (obtained for $p=q=1$):}
%
%Setting \revised{$p=1$} and $q=1$ leads to the condition
\[
	g_{k+1}^\T d_{k+1} \le -\sigma \|g_{k+1}\|^2 
	\quad \mbox{and} \quad 
	\|d_{k+1}\| \le \kappa \|g_{k+1}\|.
\]
Such a condition is typical of gradient-related directions, 
and has been instrumental in obtaining complexity guarantees for 
gradient-type methods~\cite{CCartis_PhRSampaio_PhLToint_2015}. In the 
context of nonlinear conjugate gradient method, similar properties have 
been used to establish global convergence~\cite{LGrippo_SLucidi_2005}. A case 
of particular interest to us is $d_k=-g_k$, which occurs when 
the restarting process is triggered.

%%%%%%%%%%%%%%%%%%%%%%%%%%%%%%%%%%%%%%%%%%%%%%%%%%%%%%%%%%%%%%%%%%%%%%%%%%%%%%%
\section{Complexity analysis}
\label{sec:wcc}
%%%%%%%%%%%%%%%%%%%%%%%%%%%%%%%%%%%%%%%%%%%%%%%%%%%%%%%%%%%%%%%%%%%%%%%%%%%%%%%

In this section, we derive a complexity result for our restarted variant 
of nonlinear conjugate gradient. \revised{Section~\ref{subsec:wcc:dec} 
provides the necessary assumptions as well as intermediate results, while 
Section~\ref{subsec:wcc:main} establishes and discusses complexity bounds for 
our algorithm.}

%%%%%%%%%%%%%%%%%%%%%%%%%%%%%%%%%%%%%%%%%%%%%%%%%%%%%%%%%%%%%%%%%%%%%%%%%%%%%%%
\subsection{Decrease lemmas}
\label{subsec:wcc:dec}

We make the following assumptions about the objective function of 
problem~\eqref{eq:pb}.

\begin{assumption}
\label{as:fC11L}
	The function $f$ is \revised{continuously} differentiable on $\R^n$ and its 
	gradient is $L$-Lipschitz continuous for $L>0$.
\end{assumption}

\begin{assumption}
\label{as:flow}
	There exists $\flow \in \R$ such that $f(x) \ge \flow$ for every 
	$x \in \R^n$.
\end{assumption}

\revised{Our complexity analysis relies on partitioning} the iterations 
into $\Ir \cup \It$, where
\begin{equation}
\label{eq:wcc:itpart}
	\begin{array}{lll}
		\It &= 
		&\{ k \in \N\ |\ 
		\revised{g_k^\T d_k \le -\sigma \|g_k\|^{1+p}
		\quad \mbox{and} \quad 
		\|d_k\| \le \kappa \|g_k\|^q}\} \\
		\Ir &= &\N \setminus \It.
	\end{array}
\end{equation}
\revised{For $k \ge 1$, the fact that $k \in \It$ is checked explicitly 
within our algorithm. If $k \in \Ir$, the restarting process is 
triggered and the search direction will be the negative gradient. For this 
reason, we say that $k \in \Ir$ is the index of a \emph{restarted iteration}, while 
$k \in \It$ is the index of a \emph{non-restarted iteration}.}
Depending on the nature of each iteration, we can bound the number 
of backtracking steps needed to compute a suitable stepsize. We begin by 
the non-restarted iterations, as the proof encompasses that of restarted 
iterations.

\begin{lemma}
\label{le:wcc:lst}
	Let Assumption~\ref{as:fC11L} hold, and \revised{let $k \in \It$ such 
	that $\|g_k\|> 0$.}
	Then, the line-search process terminates after at most 
	\revised{$\lfloor \bar{j}_{\It,k}+1 \rfloor$} 
	iterations, where
	\begin{equation}
	\label{eq:itslt}
		\bar{j}_{\It\revised{,k}} := \left[
		\log_{\theta}\left(\frac{2(1-\eta)\sigma}{\kappa^2 L}\right)
		\|g_k\|^{1+p-2q}\right]_+.
	\end{equation}
	Moreover, the resulting decrease at the $k$th iteration satisfies
	\begin{equation}
	\label{eq:declt}
		f(x_k) - f(x_{k+1}) > c_{\It} \min\left\{\|g_k\|^{1+p},
		\|g_k\|^{2(1+p-q)} \right\},
	\end{equation}
	where 
	\[
		c_{\It}:=\eta\sigma\min\left\{1,
		\frac{2(1-\eta)\sigma\theta}{\kappa^2 L}\right\}.
	\]
\end{lemma}

\begin{proof}
	If the decrease condition~\eqref{eq:lscond} holds for $\alpha_k=1$, 
	then the bound~\eqref{eq:itslt} holds. Moreover, 
	combining~\eqref{eq:lscond} with 
	\revised{the definition of $\It$ in~\eqref{eq:wcc:itpart}}
%	the condition~\eqref{eq:conddk}
	gives
	\begin{equation}
	\label{eq:lsggunitstep}
		f(x_k) - f(x_{k+1}) > -\eta\alpha_k g_k^\T d_k 
		\ge \eta\sigma\|g_k\|^{1+p},
	\end{equation}
	hence~\eqref{eq:declt} also holds.
	
	Suppose now that the line-search condition~\eqref{eq:lscond} fails 
	for some $\alpha = \theta^j$ with $j \in \N$. Using a Taylor expansion of 
	$f$ at $x_k$ \revised{(see, e.g.~\cite[(4.1.2)]{SJWright_2018}), we have
	\[
		f(x_k+\alpha d_k) \le f(x_k)+\alpha g_k^\T d_k + \frac{L}{2}\alpha^2 \|d_k\|^2 .
	\]
	Thus,}
	\begin{eqnarray*}
		\eta\alpha g_k^\T d_k &\le &\revised{f(x_k+\alpha d_k)}-f(x_k) \\
		&\le &\alpha g_k^\T d_k + \frac{L}{2}\alpha^2 \|d_k\|^2 \\
		&\le &\alpha g_k^\T d_k + \frac{\kappa^2 L}{2}\alpha^2 \|g_k\|^{2q},
	\end{eqnarray*}
	where the last inequality comes \revised{from~\eqref{eq:wcc:itpart}}. 
	Re-arranging the terms, we obtain
	\begin{eqnarray}
	\label{eq:contradalpha}
		-(1-\eta)\alpha g_k^\T d_k &\le &\frac{\kappa^2 L}{2}\alpha^2 \|g_k\|^{2q}
		\nonumber \\
		(1-\eta)\alpha \sigma \|g_k\|^{1+p} &\le &\frac{\kappa^2 L}{2}\alpha^2 \|g_k\|^{2q} 
		\nonumber \\
		\frac{2(1-\eta)\sigma}{\kappa^2 L}\|g_k\|^{1+p-2q} &\le &\alpha.
	\end{eqnarray}
	The condition~\eqref{eq:contradalpha} can only hold for $j \le \bar{j}_{\It\revised{,k}}$: 
	as a result, the line-search 
	process must terminate after 
	\revised{$j_k \le \lfloor \bar{j}_{\It,k}+1 \rfloor$} iterations. Moreover, 
	since the line search did not terminate after $j_k-1$ iterations, we have
	\[
		\theta^{j_k-1} \ge \frac{2(1-\eta)\sigma}{\kappa^2 L} 
		\revised{\|g_k\|^{1+p-2q}}
		\quad \Leftrightarrow \quad 
		\alpha_k = \theta^{j_k} \ge \frac{2(1-\eta)\theta\sigma}{\kappa^2 L}
		\|g_k\|^{1+p-2q}.
	\]
	Consequently, the function decrease at iteration $k$ satisfies:
	\begin{equation}
	\label{eq:funcdecnonunit}
		f(x_k)-f(x_{k+1}) > \revised{-}\eta \alpha_k g_k^\T d_k 
		\ge \frac{2\eta(1-\eta)\theta\sigma^2}{\kappa^2 L}\|g_k\|^{2(1+p-q)}
		\ge c_{\It} \|g_k\|^{2(1+p-q)},
	\end{equation}
	hence~\eqref{eq:declt} also holds in this case.
\end{proof}

We now consider the restarted iterations. In that case, we have $d_k=-g_k$, 
implying that
\begin{equation}
\label{eq:wcc:condIr}
	g_k^\T d_k = -\|g_k\|^2 \quad \mbox{and} \quad 
	\|d_k\| = \|g_k\|.
\end{equation}
Thus, for the restarted iterations, the search direction satisfies a 
property analogous to \revised{that defining non-restarted iterations 
in~\eqref{eq:wcc:itpart}} with $\sigma=\kappa=1$ and $p=q=1$. 
A reasoning identical to that used in the proof of Lemma~\ref{le:wcc:lst} 
leads to the following result.

\begin{lemma}
\label{le:wcc:lsr}
	Let Assumption~\ref{as:fC11L} hold, and \revised{let $k \in \Ir$ such 
	that $\|g_k\| > 0$.}
	Then, the line-search process terminates after at most $\bar{j}_{\Ir}+1$ 
	iterations, where
	\begin{equation}
	\label{eq:itslr}
		\bar{j}_{\Ir} := \left[
		\log_{\theta}\left(\frac{2(1-\eta)}{L}\right)
		\right]_+.
	\end{equation}
	Moreover, the resulting decrease at the $k$th iteration satisfies
	\begin{equation}
	\label{eq:declr}
		f(x_k) - f(x_{k+1}) > c_{\Ir} \|g_k\|^2,
	\end{equation}
	where 
	\[
		c_{\Ir}:=\eta\min\left\{1,
		\frac{2(1-\eta)\theta}{L}\right\}.
	\]
\end{lemma}

The result of Lemmas~\ref{le:wcc:lst} and~\ref{le:wcc:lsr} are instrumental 
to bounding the number of iterations necessary to reach an approximate 
stationary point. 

%%%%%%%%%%%%%%%%%%%%%%%%%%%%%%%%%%%%%%%%%%%%%%%%%%%%%%%%%%%%%%%%%%%%%%%%%%%%%%%
\subsection{Main results}
\label{subsec:wcc:main}

Our main result is a bound on the number of iterations performed by the 
algorithm prior to reaching an $\epsilon$-stationary point: this bound also 
applies to the number of gradient evaluations.

\begin{theorem}
\label{th:wcc:gdits}
	Let Assumptions~\ref{as:fC11L} and~\ref{as:flow} hold. 
	Suppose that $1+p-q \ge 0$. Then, the number of 
	iterations (and objective gradient evaluations) required by Algorithm~\ref{alg:algo} to reach a point 
	satisfying~\eqref{eq:epspoint} is at most
	\begin{equation}
	\label{eq:wcc:gdits}
		K_{\epsilon}:= \left\lfloor \frac{f(x_0)-\flow}{c_{\Ir}} 
		\epsilon^{-2} + \frac{f(x_0)-\flow}{c_{\It}}\epsilon^{-\max\{1+p,2(1+p-q)\}} 
		\right\rfloor.
	\end{equation}
\end{theorem}

\begin{proof}
	Let $K \in N$ be such that $\|\nabla f(x_k)\| > \epsilon$ for 
	any $k=0,\dots,K-1$. Following our partitioning~\eqref{eq:wcc:itpart}, 
	we define the index sets 
	\begin{eqnarray*}
		\It_K &:= &\It \cap \{0,\dots,K-1\},\\
		\Ir_K &:= &\Ir \cap \{0,\dots,K-1\}.
	\end{eqnarray*}
	
	For any $k \in \It_K$, the result of
	Lemma~\ref{le:wcc:lst} applies, and we have
	\begin{equation}
	\label{eq:dect}
		f(x_k) - f(x_{k+1}) \ge c_{\It} \min\left\{\|g_k\|^{1+p},
		\|g_k\|^{2(1+p-q)} \right\} \ge c_{\It} \epsilon^{\max\{1+p,2(1+p-q)\}}.
	\end{equation}
	On the other hand, if $k \in \Ir_K$, applying Lemma~\ref{le:wcc:lsr} 
	gives
	\begin{equation}
	\label{eq:decr}
		f(x_k)-f(x_{k+1}) \ge c_{\Ir} \|g_k\|^2 \ge c_{\Ir} \epsilon^2.
	\end{equation}
	
	We now consider the sum of function changes over all 
	$k \in \{0,\dots,K-1\}$. By Assumption~\ref{as:flow}, we obtain
	\begin{eqnarray*}
		f(x_0) - \flow &\ge &f(x_0) - f(x_K) \\
		&\ge &\sum_{k=0}^{K-1} \revised{\left[ f(x_k) - f(x_{k+1})\right]} \\
		&\ge &\sum_{k \in \It_K}  \revised{\left[ f(x_k) - f(x_{k+1})\right]}
		+ \sum_{k \in \Ir_K}  \revised{\left[ f(x_k) - f(x_{k+1})\right]} \\
		&> &\sum_{k \in \It_K} c_{\It} 
		\epsilon^{\max\{1+p,2(1+p-q)\}} 
		+ \sum_{k \in \Ir_K} c_{\Ir} \epsilon^2.
	\end{eqnarray*}
	Since the right-hand side consists in two sums of positive terms, the 
	above inequality implies that
	\[
		f(x_0) - \flow > \sum_{k \in \It_K} c_{\It}
		\epsilon^{\max\{1+p,2(1+p-q)\}}  
		\; \Leftrightarrow \; 
		|\It_K| < \frac{f(x_0)-\flow}{c_{\It}} \epsilon^{-\max\{1+p,2(1+p-q)\}}
	\]
	and
	\[
		f(x_0) - \flow > \sum_{k \in \Ir_K} c_{\Ir} \epsilon^2 
		\; \Leftrightarrow \;
		|\Ir_K| < \frac{f(x_0)-\flow}{c_{\Ir}} \epsilon^{-2}.
	\]
	Using $|\It_K|+|\Ir_K|=K$ finally yields 
	\[
		K < \frac{f(x_0)-\flow}{c_{\Ir}} 
		\epsilon^{-2} + \frac{f(x_0)-\flow}{c_{\It}}\epsilon^{-\max\{1+p,2(1+p-q)\}},
	\]	
	hence $K \le K_{\epsilon}$.
\end{proof}

Note that the complexity bound of Theorem~\ref{th:wcc:gdits} also guarantees 
global convergence of the algorithmic framework, since it holds for 
$\epsilon$ arbitrarily close to $0$. More precisely, it is possible to show 
that
\[
	\liminf_{k \rightarrow \infty} \|\nabla f(x_k)\| = 0,
\]
which is a typical \revised{convergence} result for nonlinear conjugate 
gradient using line search.

By combining the result of Theorem~\ref{th:wcc:gdits} with that of the 
Lemma~\ref{le:wcc:lst} and~\ref{le:wcc:lsr}, we can also provide an 
evaluation complexity bound of Algorithm~\ref{alg:algo}.

\begin{corollary}
\label{co:wccevals}
	Under the assumptions of Theorem~\ref{th:wcc:gdits}, suppose further 
	that $1+p-2q=0$. Then, the number of 
	function evaluations required by Algorithm~\ref{alg:algo} to reach a 
	point satisfying~\eqref{eq:epspoint} is at most 
	\begin{equation}
	\label{eq:wccevals}
		\left\lfloor\left[
		\log_{\theta}\left(\frac{2(1-\eta)\sigma}{\kappa^2 L}\right)
		\right]_+ +1\right\rfloor K_{\epsilon},
	\end{equation}
	where $K_{\epsilon}$ is defined in~\eqref{eq:wcc:gdits}.
\end{corollary}

\begin{proof}
	\revised{
	Since $1+p-2q=0$, we have 
	\[
		\forall k \in \It, \quad 
		\bar{j}_{\It,k} = \left[
		\log_{\theta}\left(\frac{2(1-\eta)\sigma}{\kappa^2 L}\right)
		\right]_+,
	\]
	hence this quantity is independent of the iteration index $k$. 
	Moreover,
	\[
		\max\left\{\left[
		\log_{\theta}\left(\frac{2(1-\eta)\sigma}{\kappa^2 L}\right)
		\right]_+,\bar{j}_{\Ir}\right\} 
		= \left[
		\log_{\theta}\left(\frac{2(1-\eta)\sigma}{\kappa^2 L}\right)
		\right]_+
	\]
	since $\kappa \ge 1$ and $\sigma \le 1$. As a result, any iteration 
	requires at most 
	\[
		\left\lfloor\left[
		\log_{\theta}\left(\frac{2(1-\eta)\sigma}{\kappa^2 L}\right)
		\right]_+ +1\right\rfloor
	\]
	function evaluations. Combining this number with the result of 
	Theorem~\ref{th:wcc:gdits} completes the proof.
	}
\end{proof}

\vspace*{1ex}
Note that the additional condition $1+p-2q=0$ can be replaced by a 
boundedness assumption on the gradient norm. Such an assumption also 
removes the need for the condition $1+p-q \ge 0$ in Theorem~\ref{th:wcc:gdits}.

%%%%%%%%%%%%%%%%%%%%%%%%%%%%%%%%%%%%%%%%%%%%%%%%%%%%%%%%%%%%%%%%%%%%%%%%%%%%%%%
\subsection{Interpretation of the complexity bounds}
\label{subsec:wcc:interp}

Both the iteration complexity bound and the evaluation complexity bound of 
Algorithm~\ref{alg:algo} are of order 
\begin{equation}
\label{eq:wccorder}
	\mathcal{O}\left( \epsilon^{-2}\right) + 
	\mathcal{O}\left( \epsilon^{-\max\{1+p,2(1+p-q)\}}\right).
\end{equation}
\revised{The proof of Theorem~\ref{th:wcc:gdits} shows that each term 
in the bound relates to a different form of iteration.}
The first part of the bound corresponds to negative gradient steps due 
to a restart in the algorithm. The other part corresponds to ``true'' 
conjugate gradient directions that satisfy condition~\eqref{eq:conddk}: 
thanks to this condition, we can certify a different decrease formula 
for these iterations. \revised{As a result, our complexity analysis interpolates 
between that of gradient descent and that of an ideal method where 
no restart would occur.} Note that the two terms in 
the maximum coincide when $1+p-2q=0$: this condition also guarantees 
a constant number of backtracking line-search iterations.

\begin{table}[h!]
\center
\begin{tabular}{l||l|l|l|l|l}
	$p$ & $1$ &$3/4$ &$1/2$ &$1/4$ &$0$ \\
	\hline
	$q=(1+p)/2$ & $1$ &$7/8$ &$3/4$ &$5/8$ &$1/2$ \\
	\hline
	Order $\epsilon^{-(1+p)}$ &$\epsilon^{-2}$ &$\epsilon^{-7/4}$ 
	&$\epsilon^{-3/2}$ &$\epsilon^{-5/4}$ &
	$\epsilon^{-1}$
\end{tabular}
\caption{\revised{Possible complexity values for different values of $p$ and 
$q$ satisfying $1+p-2q=0$ and $1+p-q \ge 0$.}}
\label{tab:wccorders}
\end{table}

Table~\ref{tab:wccorders} shows several possible choices for $p$ and 
$q$ that satisfy the requirements of our analysis. We focus on choices 
for which $1+p \le 2$, since those lead to a better complexity 
bound for the number of \revised{non-restarted} iterations. \revised{C}hoosing $p$ 
between $0$ and $1$ implies that the bound varies between $\epsilon^{-2}$ 
and $\epsilon^{-1}$\revised{, suggesting that} the overall number of 
iterations could be better than that of gradient descent\revised{. H}owever, 
regardless of the choice of $p$, there is 
no a priori guarantee that the number of non-restarted iterations 
will be significantly larger than that of restarted iterations. In the 
next section, we investigate this behavior in the context of nonconvex
regression problems.

%$\kappa^2 \sigma^{-2} \epsilon^{-2}$: 
%as such, they match the complexity of gradient-type methods in terms of 
%dependencies on $\epsilon$~\cite{CCartis_NIMGould_PhLToint_2010,
%CCartis_PhRSampaio_PhLToint_2015}. On the other hand, the ratio 
%$\tfrac{\kappa}{\sigma}$ illustrates the theoretical impact of our modified 
%restart condition: this ratio will typically be larger than $1$, implying that 
%our complexity bound will worsen as we allow our search directions that depart 
%more from the negative gradient (both in terms of angle and in terms of norm).

%

%\begin{remark}
%	Although the focus of this paper is on nonconvex problems, the decrease 
%	guarantee provided by Lemma~\ref{le:ls} is sufficient to derive guarantees 
%	on convex and strongly convex problems. More precisely, if we assume that 
%	problem~\eqref{eq:pb} has a minimum $f^*$ and that the objective function 
%	satisfies Assumption~\ref{as:fC11L}, it can be shown that Algorithm~\ref{alg:algo} 
%	will compute a point $x_k$ such that $f(x_k)-f^* \le \epsilon$ in at most 
%	$\mathcal{O}(\epsilon^{-1})$ iterations if the objective function is  
%	convex, and at most 
%	$\mathcal{O}\left(\tfrac{L}{\mu}\ln(\epsilon^{-1})\right)$ iterations if 
%	it is $\mu$-strongly convex. These results match their counterparts for 
%	gradient descent~\cite{SJWright_2018} in terms of dependencies on 
%	$\epsilon$, and similar comments to those above can be made regarding
%	the influence of the parameters $\sigma$ and $\kappa$.
%\end{remark}

%%%%%%%%%%%%%%%%%%%%%%%%%%%%%%%%%%%%%%%%%%%%%%%%%%%%%%%%%%%%%%%%%%%%%%%%%%%%%%%
\section{Numerical \revised{illustration}}
\label{sec:num}
%%%%%%%%%%%%%%%%%%%%%%%%%%%%%%%%%%%%%%%%%%%%%%%%%%%%%%%%%%%%%%%%%%%%%%%%%%%%%%%

In this section, we investigate the numerical behavior of our nonlinear 
conjugate gradient \revised{algorithm with the modified restarting condition.}
\revised{The purpose of these experiments is twofold. On one hand, we aim at 
better understanding the efficiency of nonlinear conjugate gradient on 
nonconvex regression problems, as demonstrated by Carmon et 
al.~\cite{YCarmon_JCDuchi_OHinder_ASidford_2017a}. We revisit this experiment 
in Section~\ref{subsec:num:pb}. On the other hand, we compare our modified 
nonlinear conjugate gradient algorithm with standard variants on a benchmark
of nonlinear optimization problems from the CUTEst 
collection~\cite{NIMGould_DOrban_PhLToint_2015}, so as to assess the numerical 
impact of enforcing complexity guarantees. This is the purpose of 
Section~\ref{subsec:num:cutest}. The setup for our experiments is described
 in Section~\ref{subsec:num:algos}.}
 
%TO DO: Interpretation of $p<0.5$?

%%%%%%%%%%%%%%%%%%%%%%%%%%%%%%%%%%%%%%%%%%%%%%%%%%%%%%%%%%%%%%%%%%%%%%%%%%%%%%%
\subsection{\revised{Algorithms and implementation}}
\label{subsec:num:algos}

\revised{Our experiments included the following methods:}
%We consider several variants on Algorithm~\ref{alg:algo} using different 
%values of the parameters $(p,q,\sigma,\kappa)$ appearing in our restart 
%condition~\eqref{eq:conddk}. We allowed the values $\sigma=0$ and $\kappa=\infty$ 
%since those correspond to a standard nonlinear conjugate gradient method 
%with a classical restart strategy.
%
%We ran our nonlinear CG methods together with the two following techniques:
\begin{itemize}
%	\item \textsc{Fixed step GD}, a version of gradient descent 
%	using a fixed stepsize $\alpha_k = \frac{1}{L}$, where $L$ is 
%	a Lipschitz constant for the problem gradient.
	\item \textsc{Semi-adaptive GD}, a version of gradient descent using an 
	adaptive estimate of the Lipschitz constant, that \revised{proceeds similarly} 
	to performing an Armijo linesearch~\cite{YCarmon_JCDuchi_OHinder_ASidford_2017a}.
	\item \textsc{Armijo GD}, a version of gradient descent using 
	\revised{the Armijo line search procedure described in Algorithm~\ref{alg:algo}}.
%	Comparing semi-adaptive GD and Armijo GD allows to 
%	assess the impact of this initialization technique.
	\item \revised{\textsc{Standard NCG}, a nonlinear conjugate gradient 
	method with Armijo line search. In this variant, restart occurs whenever 
	$g_{k+1}^\T d_{k+1} \ge 0$, in which case $d_{k+1} \leftarrow -g_{k+1}$. (Note 
	that this corresponds to setting $\sigma=0$ and $\kappa = \infty$ in  
	Algorithm~\ref{alg:algo}, but we single out this variant for simplicity.)}
	\item \revised{\textsc{Orthog NCG}, a nonlinear conjugate gradient 
	method with Armijo line search. Restarting occurs there whenever
	$|g_k^\T g_{k+1}| \ge \sigma \|g_k\|^2$ with $\sigma=0.01$. This condition 
	measures the loss of orthogonality between successive 
	gradients~\cite[Chapter 5]{JNocedal_SJWright_2006}. Although more elaborate 
	variants based on this condition can be considered~\cite{SKarimi_SAVavasis_2018}, 
	we use one that merely restarts based on a single test, akin to our proposed 
	method.}
	\item \revised{\textsc{Restarted NCG ($p$)}, our implementation of 
	Algorithm~\ref{alg:algo} with $\sigma=0.01$, $\kappa=100$, $q=\tfrac{1+p}{2}$ 
	and a variable value for $p$.}
\end{itemize}
\revised{
All variants of NCG were tested with the four formulas for $\beta_{k+1}$ given 
in~\eqref{eq:ncgFR}--\eqref{eq:ncgHZ}.
The parameters of the Armijo line search were set as $\eta=\theta=0.5$ for all 
methods. The line search was used with initial value $1$ at iteration $0$, then 
the value $2\alpha_k$ was used as initial value for iteration $k+1$. This procedure 
was not used for the Semi-adaptive GD method.
}

All the algorithms were implemented in \revised{MATLAB R2021a.}
%Python 3.7.4 using NumPy version 1.17.2, 
%SciPy version 1.3.1 and Matplotlib version 3.1.1. 
Experiments were run on 
a Dell Latitude 7400 running Ubuntu 20.04 with Intel(R) Core(TM) i7-8665U 
CPU @ 1.90GHz and 31.2GB of memory.

%%%%%%%%%%%%%%%%%%%%%%%%%%%%%%%%%%%%%%%%%%%%%%%%%%%%%%%%%%%%%%%%%%%%%%%%%%%%%%%
\subsection{\revised{Nonconvex regression tasks}}
\label{subsec:num:pb}
%%%%%%%%%%%%%%%%%%%%%%%%%%%%%%%%%%%%%%%%%%%%%%%%%%%%%%%%%%%%%%%%%%%%%%%%%%%%%%%

\revised{Our first set of experiments follows the setup of Carmon et 
al.~\cite{YCarmon_JCDuchi_OHinder_ASidford_2017a}, in which a basic nonlinear 
conjugate gradient was found to be quite efficient on randomly generated 
nonconvex regression problems.} 
Each problem instance corresponds to a dataset $\{(a_i,b_i)\}_{i=1}^m$, 
where $a_i \in \R^n$ and $b_i \in \R^m$, 
with $n=30$ and $m=60$. Every 
vector $a_i$ is generated according to a Gaussian distribution of zero mean 
and identity covariance matrix, which we denote by 
$a_i \sim \mathcal{N}(0,I_n)$, where $I_n$ is the identity matrix in 
$\R^{n \times n}$.
Letting $A = [a_i^\T] \in \R^{m \times n}$ and $b=[b_i]_{i=1}^m$, we generate 
the vector $b$ according to the following formula: 
$$
    b = A z + 3 \nu_1 + \nu_2,
$$
where $z \sim \mathcal{N}(0,4 I_n)$, $\nu_1 \sim \mathcal{N}(0,I_m)$, and 
$\nu_2 \in \mathbb{R}^m$ has i.i.d. components drawn following a Bernoulli 
distribution of parameter 0.3.

\begin{figure}
  \centering
  \includegraphics[width=0.75\textwidth]{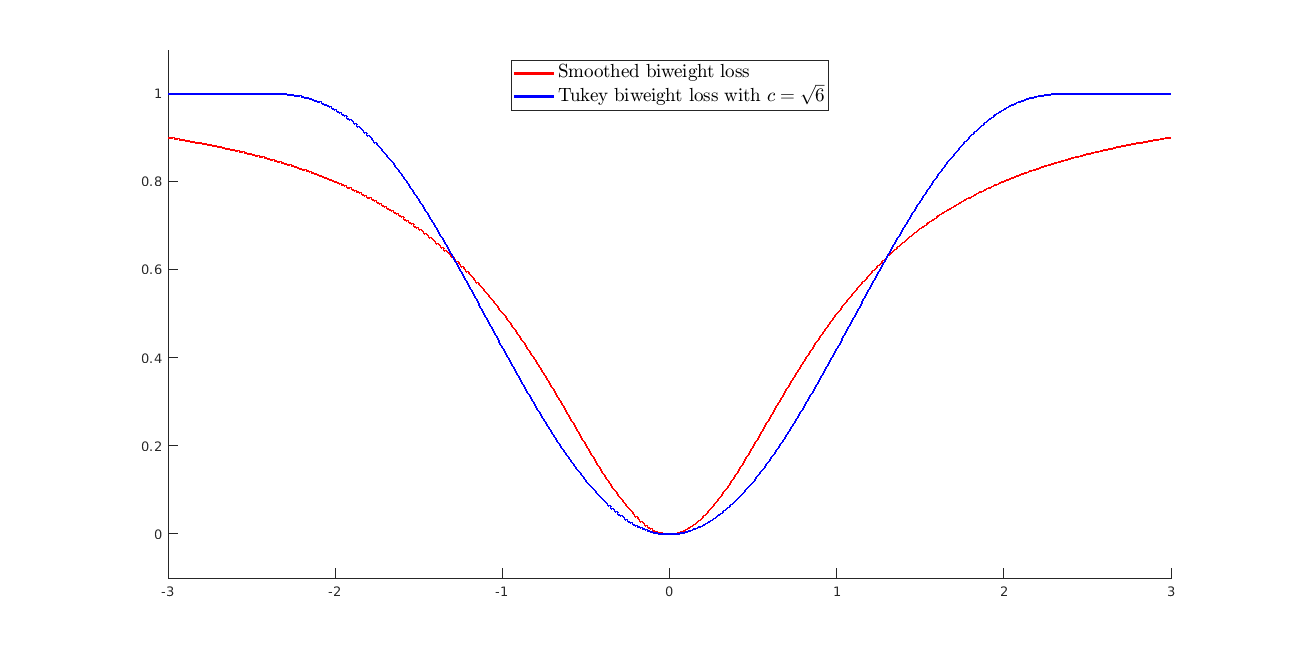}
  \caption{Nonconvex losses for robust linear regression.}
  \label{fig:losses}
\end{figure}

Given a dataset $\{(a_i,b_i)\}_{i=1}^m$ of this form, we consider the robust 
regression problem:
\begin{equation}
\label{eq:BLpb}
	\min_{x \in \R^d} f_{SB}(x):=\frac{1}{m}\sum_{i=1}^m \phi(a_i^\T x - b_i),
\end{equation}
where $\phi:\R \rightarrow [0,\infty)$ is a smoothed biweight loss 
function defined by 
\begin{equation*}
%\label{eq:biweightloss}
	\phi(t) = \frac{t^2}{1+t^2}.
\end{equation*}
Carmon et al~\cite{YCarmon_JCDuchi_OHinder_ASidford_2017a} introduced function 
$\phi$ as a smooth proxy for the Tukey biweight loss function~\cite{AEBeaton_JWTukey_1974}. 
This terminology denotes a family of functions parameterized by $c>0$ as follows: 
\[
	\forall t \in \R, \qquad
	\rho_c(t) = 
	\left\{
		\begin{array}{ll}
			\tfrac{t^6}{6 c^4}-\tfrac{t^4}{2 c^2}+\tfrac{t^2}{2} &\mbox{if\ } |t| \le c, \\
			 & \\
			\frac{c^2}{6} &\mbox{otherwise.}
		\end{array}
	\right.
\]  
Given the same data $\{(a_i,b_i)\}_{i=1}^m$ than that used to 
define problem~\eqref{eq:BLpb}, we also consider the problem
\begin{equation}
\label{eq:TBpb}
	\min_{x \in \R^d} f_{TB}(x):=\frac{1}{m}\sum_{i=1}^m \rho_{\sqrt{6}}(a_i^\T x - b_i),
\end{equation}
\revised{where we chose} $c=\sqrt{6}$ so as to be close to the smoothed version of the 
biweight loss.
The solution of problem~\eqref{eq:TBpb} belongs to the class of robust 
M-estimators in statistics; compared to a standard linear least-squares 
formulation, problem~\eqref{eq:TBpb} is more robust to outliers in the data.
Figure~\ref{fig:losses} shows the shape of $\phi$ and 
$\rho_{\sqrt{6}}$. Both losses are nonconvex, resulting in both functions $f_{SB}$ and $f_{TB}$ being nonconvex. Those functions also satisfy Assumption~\ref{as:fC11L}.
%, respectively 
%with Lipschitz constants $2\max_{i=1,\dots,m} \|a_i\|^2$ and 
%$12\max_{i=1,\dots,m} \|a_i\|^2$. 
\revised{Note that both functions are twice continuously differentiable, but 
that only $f_{SB}$ is infinitely smooth.}

%For every $x \in \mathbb{R}^d$, we have
%\begin{equation}
%\label{eq:pbgrad}
%    \nabla f(x) = \frac{2}{m} \sum_{i=1}^m \frac{a_i^T x -b_i}
%    {(1+(a_i^T x - b)^2)^2}\,a_i = \frac{2}{m} A^T c \in \R^d,
%\end{equation}
%where 
%\[
%	A = [a_i^\T] \in \R^{m \times d} \quad \mbox{and} \quad 
%	c= \left[ \frac{a_i^T x -b_i}
%    {(1+(a_i^T x - b_i)^2)^2} \right] \in \R^m
%\]
%as well as
%\[
%    \nabla^2 f(x) = \frac{2}{m} A^T \mathrm{diag}\left(
%    \left\{\frac{1-3(a_i^T x - b_i)^2}{(1+(a_i^T x - b_i)^2)^3}
%    \right\}_{i=1}^m \right) A \in \R^{d \times d}.
%\]
%where $\mathrm{diag}(v)$ is the diagonal matrix whose diagonal entries are 
%the coefficients of $v$.
%
%\begin{remark}
%For every $x \in \mathbb{R}^d$, we have:
%\[
%	\|\nabla^2 f(x)\| \le \frac{2}{m}\|A^T\|\|A\|, 
%\]
%therefore the function has a Lipschitz continuous gradient.
%\end{remark}
%
%
%\begin{example}
%\label{ex:pbtoyex}
%	Let $m=4$ and $d=2$. Given $\beta \in \R$, we consider the dataset
%	\[
%		A = \begin{bmatrix}
%				1 &0 \\
%				0 &1 \\
%				-1 &0 \\
%				0 &-1
%			\end{bmatrix}, 
%		\qquad
%		b = \begin{bmatrix}
%			1 \\
%			1 \\
%			\beta \\
%			-1
%		\end{bmatrix}.
%	\]
%	When $\beta=-1$, the problem has one global optimum $x^* = [1\ 1]^\T$. 
%	When $\beta \neq -1$, the problem possesses one saddle point at 
%	$\left[ \tfrac{1-\beta}{2}\ 1 \right]^\T$.
%\end{example}

\paragraph{\revised{Results}} 
\revised{We generated 1000 datasets for linear regression according to the 
procedure above. For each dataset, we consider the associated nonconvex 
regression problems~\eqref{eq:BLpb} and~\eqref{eq:TBpb}, and compare gradient 
descent techniques with standard NCG and several instances of restarted NCG 
(Algorithm~\ref{alg:algo}).
For every problem, all methods were run until $\|g_k\| \le \epsilon=10^{-4}$ or 
a budget of 10000 iterations was exhausted.
 We conducted experiments with values of $p$ uniformly 
distributed between $0$ and $1$. In all cases, the performance was consistently 
better for $p \ge 0.5$, with best behavior typically obtained for $p \in [0.5,0.75]$. 
We thus elected to present results using $p \in \{0,0.25,0.5,0.75,1\}$, as those 
values are representative of the behavior of our framework.}

\begin{figure}
  \centering
  \includegraphics[width=\textwidth]{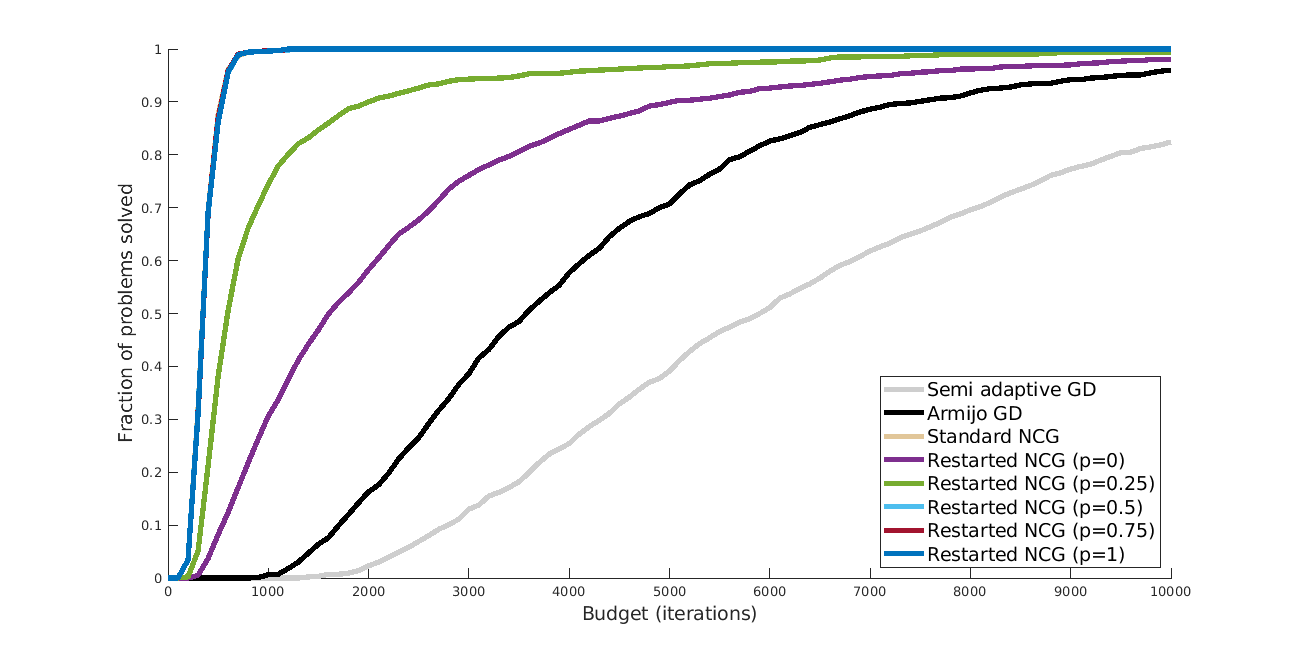}
  \caption{\revised{Fraction of nonconvex problems with the smooth biweight 
  loss~\eqref{eq:BLpb} solved as function of the iteration budget. 
  The curves corresponding to Restarted NCG 
  with $p \in \{0.5,0.75,1\}$ 
  overlap with that of Standard NCG.
  All NCG variants use the PRP+ formula~\eqref{eq:ncgPRP+}.}}
  \label{fig:dataprofBL:PRPplus}
\end{figure}

\revised{We begin by considering the original setup of Carmon et al.~\cite{YCarmon_JCDuchi_OHinder_ASidford_2017a}, where nonlinear CG was applied 
to problem~\eqref{eq:BLpb} using the PRP+ formula~\eqref{eq:ncgPRP+}. 
Table~\ref{tab:statsBL:PRPplus} shows the average percentage 
of restarted iterations among all iterations for each of the tested methods. The standard 
NCG method has a remarkably low percentage of restarted iterations, i.e. a very small 
fraction of iterations produced directions that were not of descent type. When we use the 
modified restarting condition of Algorithm~\ref{alg:algo} with $p<0.5$, we observe that 
restarted iterations form the majority of iterations on average. This behavior suggests 
that the algorithm tends to use gradient descent directions, and that most directions 
produced by a nonlinear CG update satisfy the restart condition~\eqref{eq:restartcond} 
when $p < 0.5$. On the contrary, when $p \ge 0.5$, the percentage of restarted iterations 
decreases significantly, indicating that the method relies on directions close to that 
of a standard nonlinear conjugate gradient method. Figure~\ref{fig:dataprofBL:PRPplus} 
shows the fraction of instances solved as a function of the iteration budget, under the 
form of data profiles~\cite{JJMore_SMWild_2009}. This figure confirms that the variants of Restarted NCG for $p \ge 0.5$ exhibit a behavior close to that of Standard NCG.}

\begin{table}[h!]
\centering
\begin{tabular}{|c|ll|}
\hline
Method &Problems solved &Avg. restart (\%) \\
\hline 
Standard NCG &1000 &0.74\% \\
\hline
\hline
NCG($0$) &1000 &83.5\% \\
NCG($0.25$) &1000 &53.2\% \\
NCG($0.5$) &1000 &0.89\% \\
NCG($0.75$) &1000 &0.76\% \\
NCG($1$) &1000 &0.76\% \\
\hline
\end{tabular}
\caption{Statistics on 1000 random instances of robust linear regression using smoothed biweight loss (problem~\eqref{eq:BLpb}) and a budget of 10000 iterations. 
\revised{All variants use the PRP+ formula~\eqref{eq:ncgPRP+}}.}
\label{tab:statsBL:PRPplus}
\end{table}

\begin{figure}
  \centering
  \includegraphics[width=1.1\textwidth]{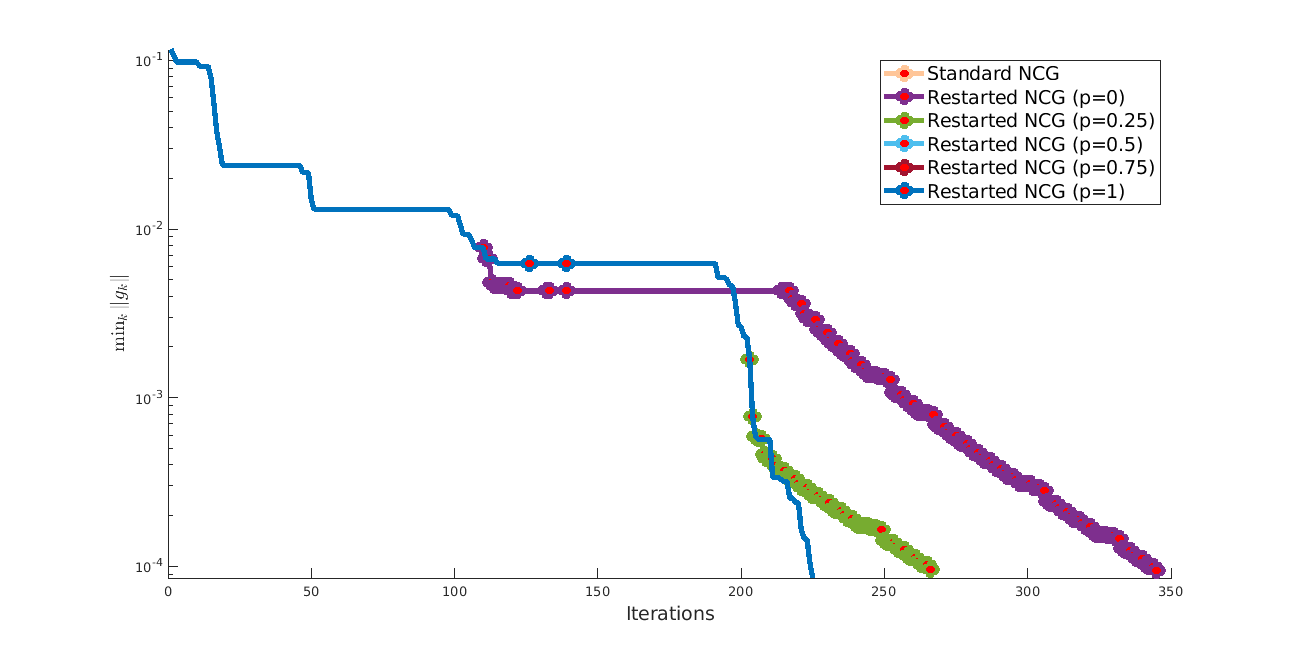}
  \caption{\revised{Minimum gradient norm for a representative run using the 
  smooth biweight loss~\eqref{eq:BLpb} using the PRP+ formula~\eqref{eq:ncgPRP+}.
  The curves corresponding to Restarted NCG 
  with $p \in \{0.5,0.75,1\}$ 
  overlap with that of Standard NCG, and all have the same 2 restarted iterations (blue 
  circles with red filling). Restarted NCG($p=0.25$) variant had 66 restarted iterations, while 
  Restarted NCG($p=0$) had 148 restarted iterations.}}
  \label{fig:dataprofBL:mingradPRPplus414}
\end{figure}

\revised{Figure~\ref{fig:dataprofBL:mingradPRPplus414} illustrates 
the behavior of the minimum gradient norm on a representative run. The number 
and location of the restarted iterations (red circles) is the same for 
Standard NCG and the restarted NCG variants for $p \ge 0.5$, suggesting that the restarting 
condition~\eqref{eq:restartcond} tends to be triggered for non-descent directions. On 
the contrary, using $p \in \{0,0.25\}$ leads to series of restarted iterations, from 
which the algorithm does not appear to recover, in that the restarting condition 
keeps being triggered. These results strongly suggest that 
most directions produced by nonlinear CG do not satisfy condition~\eqref{eq:conddk} when 
$p$ is small. However, we observe that in the early iterations, all methods perform 
non-restarted iterations, hence our conditions appear to be satisfied at the beginning 
of every run.}

\begin{figure}
  \centering
  \includegraphics[width=0.8\textwidth]{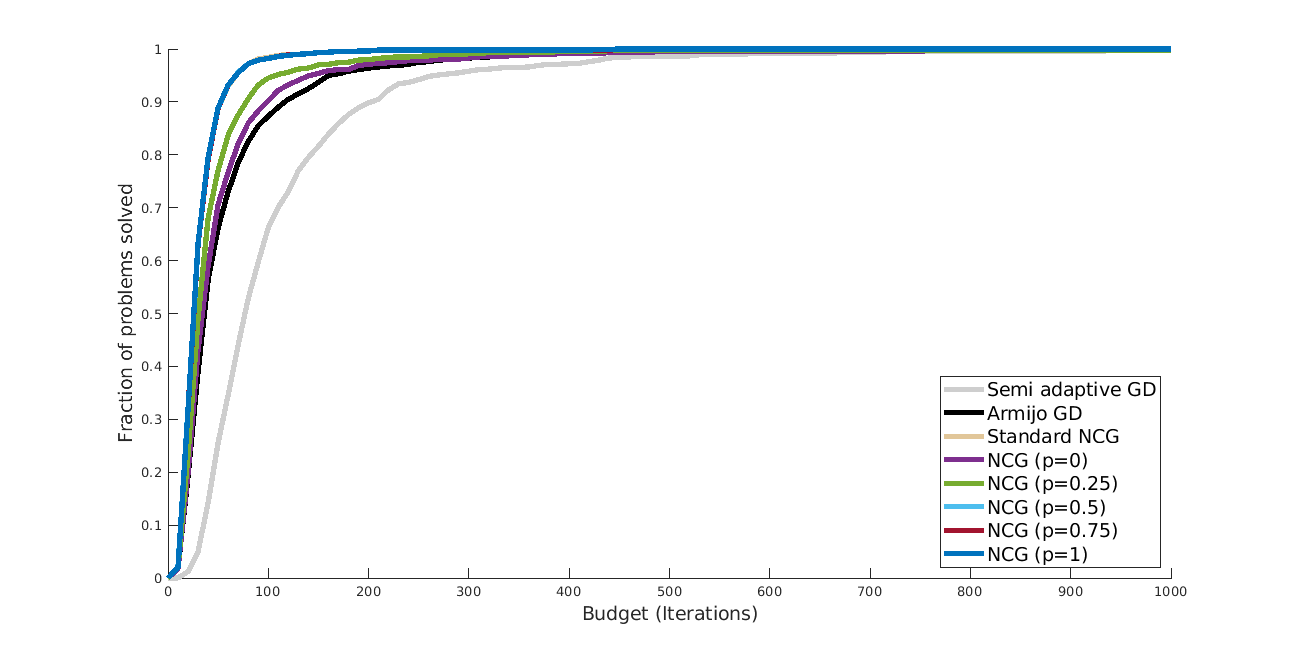}
  \caption{
  \revised{Fraction of nonconvex problems with the Tukey biweight 
  loss~\eqref{eq:TBpb} solved as function of the iteration budget (truncated 
  at 1000). 
  The curves corresponding to Restarted NCG 
  with $p \in \{0.5,0.75,1\}$ 
  overlap with that of Standard NCG.
  All NCG variants use the PRP+ formula~\eqref{eq:ncgPRP+}.}}
  \label{fig:dataprofTB:PRPplus}
\end{figure}

\begin{figure}
  \centering
  \includegraphics[width=1.1\textwidth]{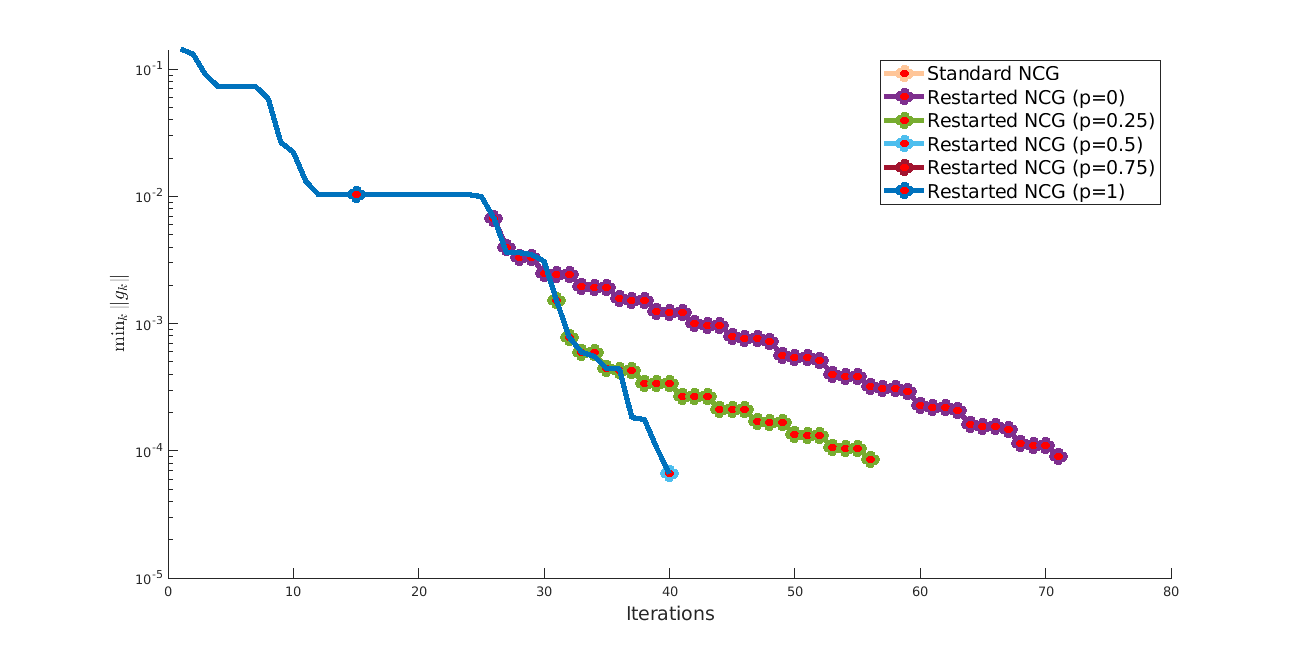}
  \caption{\revised{Minimum gradient norm for a representative run using the 
  smooth biweight loss~\eqref{eq:BLpb} using the PRP+ formula~\eqref{eq:ncgPRP+}.
  The curves corresponding to Restarted NCG 
  with $p \in \{0.75,1\}$ 
  overlap with that of Standard NCG, and all have the same restarted iteration (red 
  circles). Restarted NCG($p=0.5$) had one extra restarted iteration, but overlaps with 
  Standard NCG otherwise. Restarted NCG($p=0.25$) variant had 27 restarted iterations, while 
  Restarted NCG($p=0$) had 47 restarted iterations.}}
  \label{fig:dataprofTB:mingradPRPplus414}
\end{figure}

\revised{Table~\ref{tab:statsTT:PRPplus} as well as 
Figures~\ref{fig:dataprofTB:PRPplus} and~\ref{fig:dataprofTB:mingradPRPplus414}
are the counterparts of the previous results for the 
Tukey biweight loss problem~\eqref{eq:TBpb}. Although the percentages of restarted iterations are 
smaller than that of the previous table, the same overall trend can be observed, with 
a reduction in the number of restarted iterations as $p$ increases, and a sharp decrease 
for $p \ge 0.5$. Figure~\ref{fig:dataprofTB:mingradPRPplus414} also confirms that 
the restarted iterations for $p \ge 0.5$ occur at the same index as that of 
Standard NCG. We note that the discrepancy between nonlinear 
conjugate gradient methods and gradient descent methods is larger on 
problem~\eqref{eq:BLpb}, where the smoothed biweight function is used. Since 
the latter function is infinitely smooth, it is possible that high-order derivatives 
are relevant for optimization, in a way that nonlinear conjugate gradient better 
captures. We point out that the complexity analysis of 
Section~\ref{sec:wcc} only assumes that the objective function is continuously 
differentiable, but that additional smoothness can improve guarantees of existing 
schemes~\cite{CCartis_NIMGould_PhLToint_2019a,
YCarmon_JCDuchi_OHinder_ASidford_2021}.
}

\begin{table}[h!]
\centering
\begin{tabular}{|c|ll|}
\hline
Method &Problems solved &Avg. restart (\%) \\
\hline 
Standard NCG &1000 &0.58\% \\
\hline
\hline
NCG($0$) &1000 &62.7\% \\
NCG($0.25$) &1000 &44.6\% \\
NCG($0.5$) &1000 &3.47\% \\
NCG($0.75$) &1000 &0.61\% \\
NCG($1$) &1000 &0.63\% \\
\hline
\end{tabular}
\caption{Statistics on 1000 random instances of robust linear regression using Tukey loss (problem~\eqref{eq:TBpb}) and a budget of 10000 iterations. 
\revised{All methods use the PRP+ formula~\eqref{eq:ncgPRP+}}.}
\label{tab:statsTT:PRPplus}
\end{table}

\revised{Using the PR formula~\eqref{eq:ncgPR}
yields very similar results to those for the PRP+ formula~\eqref{eq:ncgPRP+}, therefore we do 
not report these results here. However, we provide results using the Hager-Zhang 
formula~\eqref{eq:ncgHZ} in Tables~\ref{tab:statsBL:HZ} and~\ref{tab:statsTT:HZ}. When this parameter 
formula is used, restart no longer occurs for Standard NCG as the direction is always 
a descent one~\cite{WWHager_HZhang_2005}. Still, for the restarted NCG variants, 
we again observe that the average percentage of restarted iterations diminishes as 
we increase the value of $p$, with $p \ge 0.5$ leading to significantly less restarts 
on average.}

%\begin{table}[h!]
%\centering
%\begin{tabular}{|c|ll|}
%\hline
%Method &Problems solved &Avg. restart (\%) \\
%\hline 
%Standard NCG &1000 &0.71\% \\
%\hline
%\hline
%NCG($0$) &973 &83.7\% \\
%NCG($0.25$) &993 &53.4\% \\
%NCG($0.5$) &1000 &1.02\% \\
%NCG($0.75$) &1000 &0.75\% \\
%NCG($1$) &1000 &0.75\% \\
%\hline
%\end{tabular}
%\caption{Statistics on 1000 random instances of robust linear regression using smoothed biweight loss (problem~\eqref{eq:BLpb}) and a budget of 10000 iterations. 
%\revised{All variants of Algorithm~\ref{alg:algo} use $\beta_{k+1}=\beta_{k+1}^{PRP}$}.}
%\label{tab:statsBL:PRP}
%\end{table}
%
%\begin{table}[h!]
%\centering
%\begin{tabular}{|c|ll|}
%\hline
%Method &Problems solved &Avg. restart (\%) \\
%\hline 
%Standard NCG &1000 &0.58\% \\
%\hline
%\hline
%NCG($0$) &1000 &62.7\% \\
%NCG($0.25$) &1000 &44.4\% \\
%NCG($0.5$) &1000 &4.66\% \\
%NCG($0.75$) &1000 &0.62\% \\
%NCG($1$) &1000 &0.61\% \\
%\hline
%\end{tabular}
%\caption{Statistics on 1000 random instances of robust linear regression using Tukey loss (problem~\eqref{eq:TBpb}) and a budget of 10000 iterations. 
%\revised{All variants of Algorithm~\ref{alg:algo} use $\beta_{k+1}=\beta_{k+1}^{PRP}$}.}
%\label{tab:statsTT:PRP}
%\end{table}
%
\begin{table}[h!]
\centering
\begin{tabular}{|c|ll|}
\hline
Method &Problems solved &Avg. restart (\%) \\
\hline 
Standard NCG &1000 &0.00\% \\
\hline
\hline
NCG($0$) &1000 &52.8\% \\
NCG($0.25$) &1000 &21.8\% \\
NCG($0.5$) &1000 &0.56\% \\
NCG($0.75$) &1000 &0.62\% \\
NCG($1$) &1000 &0.76\% \\
\hline
\end{tabular}
\caption{\revised{Statistics on 1000 random instances of robust linear regression using smoothed biweight loss (problem~\eqref{eq:BLpb}) and a budget of 10000 iterations. 
All methods use the HZ formula~\eqref{eq:ncgHZ}}.}
\label{tab:statsBL:HZ}
\end{table}

\begin{table}[h!]
\centering
\begin{tabular}{|c|ll|}
\hline
Method &Problems solved &Avg. restart (\%) \\
\hline 
Standard NCG &1000 &0.00\% \\
\hline
\hline
NCG($0$) &1000 &48.5\% \\
NCG($0.25$) &1000 &26.8\% \\
NCG($0.5$) &1000 &1.28\% \\
NCG($0.75$) &1000 &0.75\% \\
NCG($1$) &1000 &0.86\% \\
\hline
\end{tabular}
\caption{\revised{Statistics on 1000 random instances of robust linear regression using Tukey loss (problem~\eqref{eq:TBpb}) and a budget of 10000 iterations. 
All methods use the HZ formula~\eqref{eq:ncgHZ}}.}
\label{tab:statsTT:HZ}
\end{table}

\begin{table}[h!]
\centering
\begin{tabular}{|c|ll|}
\hline
Method &Problems solved &Avg. restart (\%) \\
\hline 
Standard NCG &9 &0.03\% \\
\hline
\hline
NCG($0$) &122 &2.98\% \\
NCG($0.25$) &197 &0.94\% \\
NCG($0.5$) &216 &0.02\% \\
NCG($0.75$) &368 &0.03\% \\
NCG($1$) &514 &0.03\% \\
\hline
\end{tabular}
\caption{Statistics on 1000 random instances of robust linear regression using smoothed biweight loss (problem~\eqref{eq:BLpb}) and a budget of 10000 iterations. 
\revised{All variants of Algorithm~\ref{alg:algo} use $\beta_{k+1}=\beta_{k+1}^{FR}$}.}
\label{tab:statsBL:FR}
\end{table}

\begin{table}[h!]
\centering
\begin{tabular}{|c|ll|}
\hline
Method &Problems solved &Avg. restart (\%) \\
\hline 
Standard NCG &629 &0.07\% \\
\hline
\hline
NCG($0$) &730 & 11.0\\
NCG($0.25$) &759 &4.59\% \\
NCG($0.5$) &769 &0.11\% \\
NCG($0.75$) &839 &0.06\% \\
NCG($1$) &876 &0.07 \\
\hline
\end{tabular}
\caption{Statistics on 1000 random instances of robust linear regression using Tukey loss (problem~\eqref{eq:TBpb}) and a budget of 10000 iterations. 
\revised{All variants of Algorithm~\ref{alg:algo} use $\beta_{k+1}=\beta_{k+1}^{FR}$}.}
\label{tab:statsTT:FR}
\end{table}

\begin{figure}
  \centering
  \includegraphics[width=\textwidth]{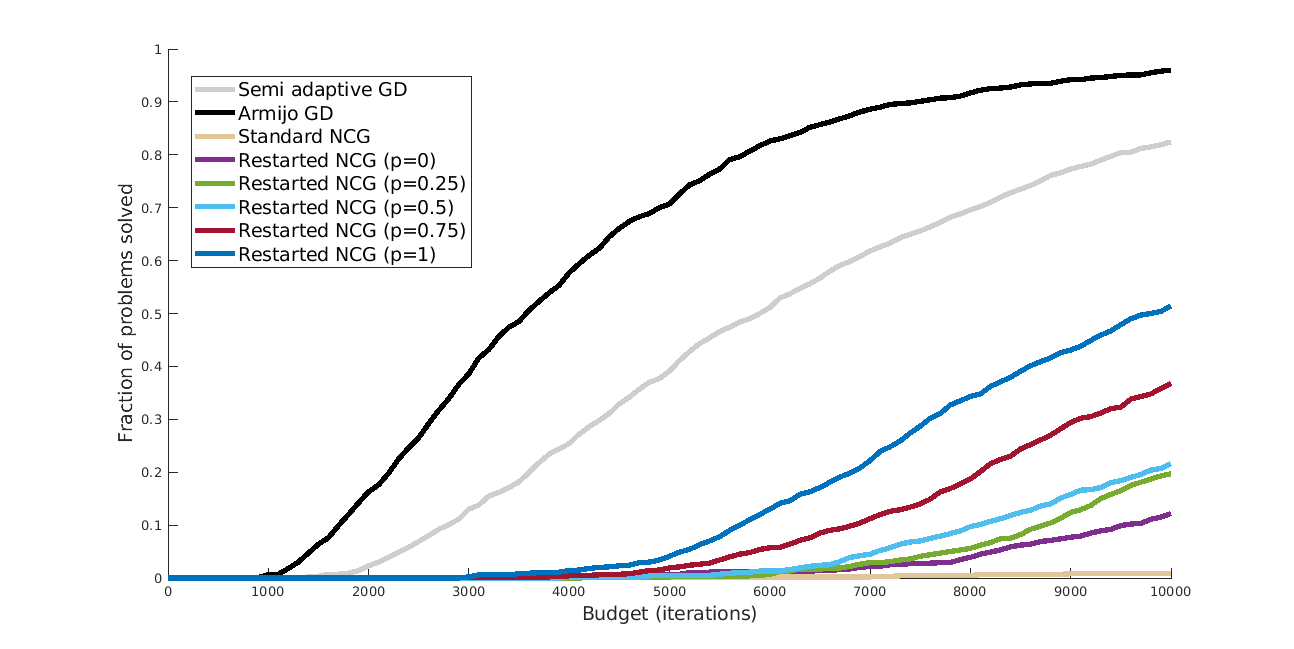}
  \caption{\revised{Fraction of nonconvex problems with the smooth biweight 
  loss~\eqref{eq:BLpb} solved as function of the iteration budget. 
  The curves corresponding to Restarted NCG 
  with $p \in \{0.5,0.75,1\}$ 
  overlap with that of Standard NCG.
  All NCG variants use the $\beta_{k+1}^{FR}$ formula.}}
  \label{fig:dataprofBL:FR}
\end{figure}

\begin{figure}
  \centering
  \includegraphics[width=\textwidth]{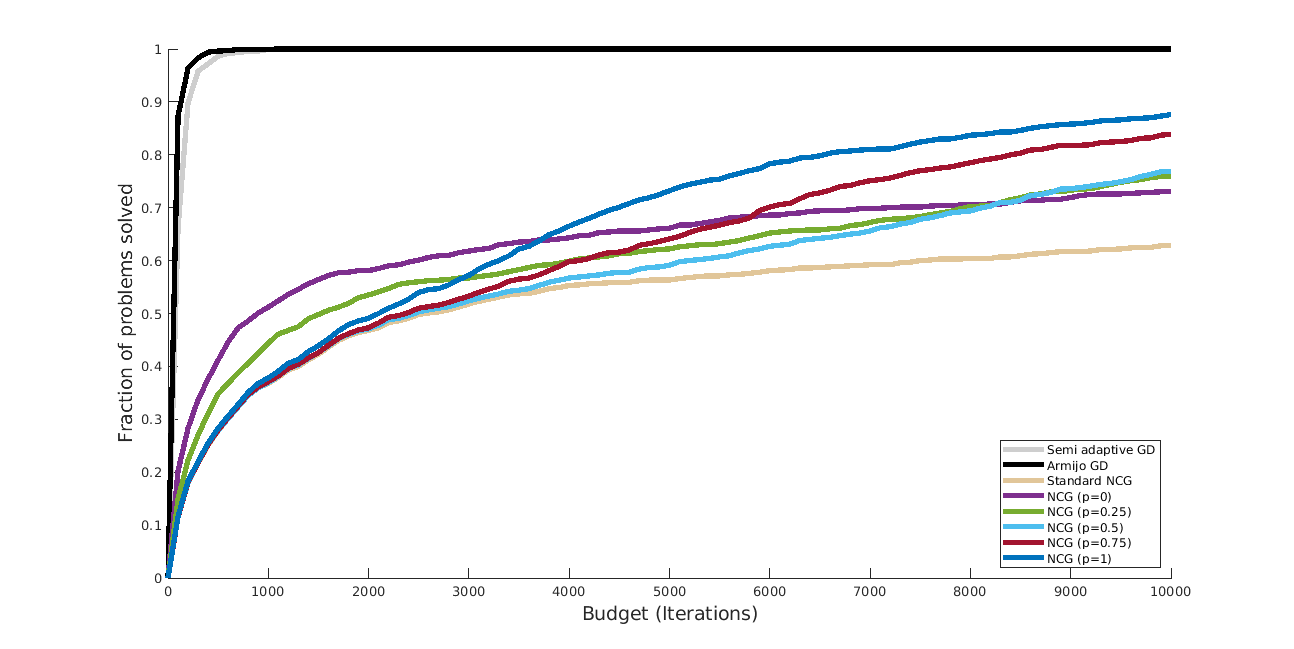}
  \caption{
  \revised{Fraction of nonconvex problems with the Tukey biweight 
  loss~\eqref{eq:TBpb} solved as function of the iteration budget. 
  The curves corresponding to Restarted NCG 
  with $p \in \{0.5,0.75,1\}$ 
  overlap with that of Standard NCG.
  All NCG variants use the $\beta_{k+1}^{FR}$ formula.}}
  \label{fig:dataprofTB:FR}
\end{figure}

\revised{
Finally, we present the results for the Fletcher-Reeves formula~\eqref{eq:ncgFR}, 
classically established as a variant with theoretical guarantees but less practical 
appeal~\cite{WWHager_HZhang_2006b}. As shown in Tables~\ref{tab:statsBL:FR} and 
\ref{tab:statsTT:FR}, this update leads to a much worse performance 
for all the NCG methods, and in that setting all methods with modified restart 
give better results than Standard NCG. Interestingly, Figures~\ref{fig:dataprofBL:FR} 
and~\ref{fig:dataprofTB:FR} show that gradient descent actually outperforms 
the NCG variants, but that adding the modified restarted condition consistently 
improves the method's performance as $p$ gets closer to 1. These results suggest 
that restarting conditions may be a way to improve the performance of a 
Fletcher-Reeves nonlinear conjugate gradient 
method.}

\revised{Overall, our experiments indicate that our proposed framework can 
perform quite closely to a standard implementation of nonlinear CG on 
nonconvex regression problems. Setting $p \ge 0.5$ emerges as the best choice to track the behavior 
of standard nonlinear CG for certain parameter formulas, and even improve 
over it in the case of the FR formula. When the PRP+ formula is used, the 
restart condition is rarely triggered, thus the number of non-restarted 
iterations dominates that of restarted iterations. Our complexity guarantees 
in this (admittedly very specific) setting can therefore improve over the theoretically 
fastest methods based on accelerated 
gradient~\cite{YCarmon_JCDuchi_OHinder_ASidford_2017a}, and provide an 
ad-hoc justification for this practical behavior.}

%%%%%%%%%%%%%%%%%%%%%%%%%%%%%%%%%%%%%%%%%%%%%%%%%%%%%%%%%%%%%%%%%%%%%%%%%%%%%%%
\subsection{\revised{Smooth optimization benchmark}}
\label{subsec:num:cutest}

\revised{
We now consider a more substantial test set consisting of 
CUTEst collection~\cite{NIMGould_DOrban_PhLToint_2015}\footnote{
\revised{Downloaded 
from GitHub on June 8, 2021.}}. Our test comprises 
149 smooth unconstrained optimization 
problems previously used to compare methods with and without complexity guarantees~\cite{FECurtis_DPRobinson_CWRoyer_SJWright_2021}. We chose the 
default dimension of each problem when indicated, or used the largest dimension 
below 1000 otherwise. The complete problem list is given in Table~\ref{tab:cutest}.}

%\revised{
\begin{tiny}
\begin{table}[h!]
\centering
\begin{tabular}{|cc|cc|cc|}
\hline
ALLINITU    &4 &ARGLINA		&200 &ARGLINB		&200 \\
ARGLINC		&200 &ARWHEAD		&1000 &BARD        &3 \\
BDEXP		&100	&BDQRTIC &100 &BEALE	&2 \\
BIGGS6      &6 &BOX3        &3 &BRATU1D		&77 \\
BRKMCC		&2 &BROWNAL     &10 &BROWNBS		&2 \\
BROWNDEN	&4
&BROYDN7D    &10
&BRYBND      &10 \\
CHAINWOO	&1000 
&CHNROSNB    &10
&CLIFF		&2 \\
CLPLATEA	&49
&CLPLATEB	&49
&CLPLATEC	&49 \\
COSINE		&1000
&CRAGGLVY	&4
&CUBE		&2 \\
CURLY10		&1000
&CURLY20		&1000
&CURLY30		&1000 \\
DECONVU		&63
&DENSCHNA    &2 
&DENSCHNB    &2 \\ 
DENSCHNC    &2 
&DENSCHND    &3 
&DENSCHNE    &3 \\
DENSCHNF    &2 
&DIXMAANA    &15
&DIXMAANB    &15 \\
DIXMAANC    &15
&DIXMAAND    &15
&DIXMAANE    &15 \\
DIXMAANF    &15
&DIXMAANG    &15
&DIXMAANH    &15 \\
DIXMAANI    &15
&DIXMAANJ    &15 
&DIXMAANK    &15 \\
DIXMAANL    &15 
&DIXON3DQ	&10 
&DJTL		&2 \\
DQDRTIC		&10 
&DQRTIC		&500
&EDENSCH		&36 \\
EG2			&1000 
&EIGENALS	&110
&EIGENBLS	&110 \\
ENGVAL1		&2 
&ENGVAL2     &3
&ERRINROS    &50 \\
EXPFIT      &2 
&EXTROSNB	&5
&FLETCBV2	&10 \\
FLETCBV3	&10 
&FLETCHBV	&10
&FLETCHCR	&10 \\
FMINSRF2	&64 
&FMINSURF    &16
&FREUROTH	&2\\
GENHUMPS	&1000 
&GENROSE		&500
&GROWTHLS    &3 \\
GULF		&3 
&HAIRY       &2  
&HATFLDD     &3 \\
HATFLDE     &3 
&HEART6LS	&6 
&HEART8LS    &8 \\
HELIX       &3 
&HILBERTA	&10
&HILBERTB	&5 \\
HIMMELBB    &2 
&HIMMELBF	&4 
&HIMMELBG    &2 \\
HIMMELBH	&2 
&HUMPS       &2
&INDEF		&1000 \\
JENSMP		&2 
&KOWOSB      &4 
&LIARWHD		&36 \\
LOGHAIRY    &2 
&MANCINO     &100
&MARATOSB    &2 \\
MEXHAT		&2 
&MEYER3	    &3
&MSQRTALS    &4  \\
MSQRTBLS    &9 
&NONCVXU2	&1000
&NONCVXUN	&1000 \\
NONDIA		&1000 
&NONDQUAR	&100
&NONMSQRT	&9 \\
OSBORNEA	&5 
&OSBORNEB	&11 
&PALMER1C	&8 \\
PALMER1D	&7 
&PALMER1E	&8
&PALMER2C	&8 \\
PALMER2E	&8 
&PALMER3C	&8
&PALMER3E	&8 \\
PALMER4C	&8 
&PALMER4E	&8
&PALMER5C	&6 \\
PALMER5D	&4 
&PALMER6C	&8
&PALMER7C	&8 \\
PALMER8C	&8 
&PENALTY1	&1000
&PENALTY2	&4 \\
PENALTY3	&50 
&PFIT1LS		&3
&PFIT2LS		&3 \\
PFIT3LS		&3 
&PFIT4LS		&3
&POWER		&10 \\
QUARTC		&25 
&ROSENBR		&2
&SCOSINE		&1000\\
SCURLY10	&1000 
&SCURLY20	&1000
&SCURLY30	&1000 \\
SINEVAL		&2 
&SINQUAD		&5
&SISSER		&2 \\
SNAIL		&2 
&SPARSINE	&1000
&SPMSRTLS	&28\\
SROSENBR	&10 
&STRATEC		&10
&TOINTQOR	&50\\
TRIDIA		&30 
&VARDIM		&10
&VAREIGVL	&50\\
VIBRBEAM	&8 
&WATSON		&31
&WOODS		&4 \\
YFITU		&3 
&ZANGWIL2	&2 
& & \\
\hline
\end{tabular}
\caption{List of CUTEst problems}
\label{tab:cutest}
\end{table}
\end{tiny}

\revised{
We only consider nonlinear conjugate gradient techniques in our comparison. 
We ran Standard NCG, Orthog NCG and Restarted NCG for $p \in \{0,0.5,0.75,1\}$ 
using different values for the $\beta_{k+1}$ parameter. All methods were run 
with a budget of 10000 iterations. A run was considered convergent whenever it 
reached a point $x_k$ such that
\[
	\|\nabla f(x_k) \| \le \epsilon \max\left\{1,\|\nabla f(x_0)\|\right\},
\]
where $\epsilon=10^{-5}$.}

\paragraph{\revised{Results}}

\begin{figure}
  \centering
  \includegraphics[width=\textwidth]{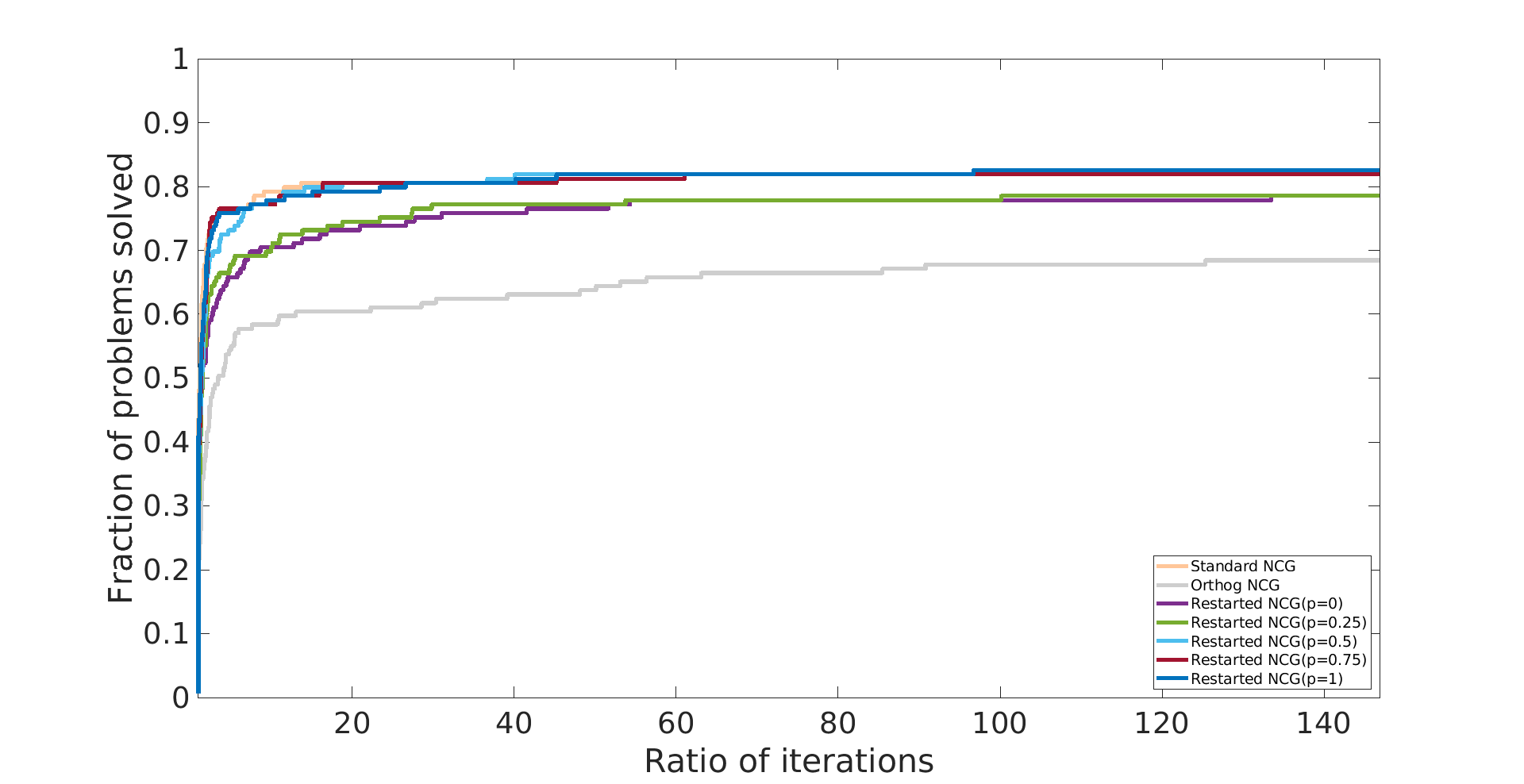}
  \caption{
  \revised{Performance profiles (iterations) 
  for nonlinear CG methods with the FR formula~\eqref{eq:ncgFR}
formula on a benchmark from CUTEst.}}
  \label{fig:cutest:FRits}
\end{figure}

\begin{figure}
  \centering
  \includegraphics[width=\textwidth]{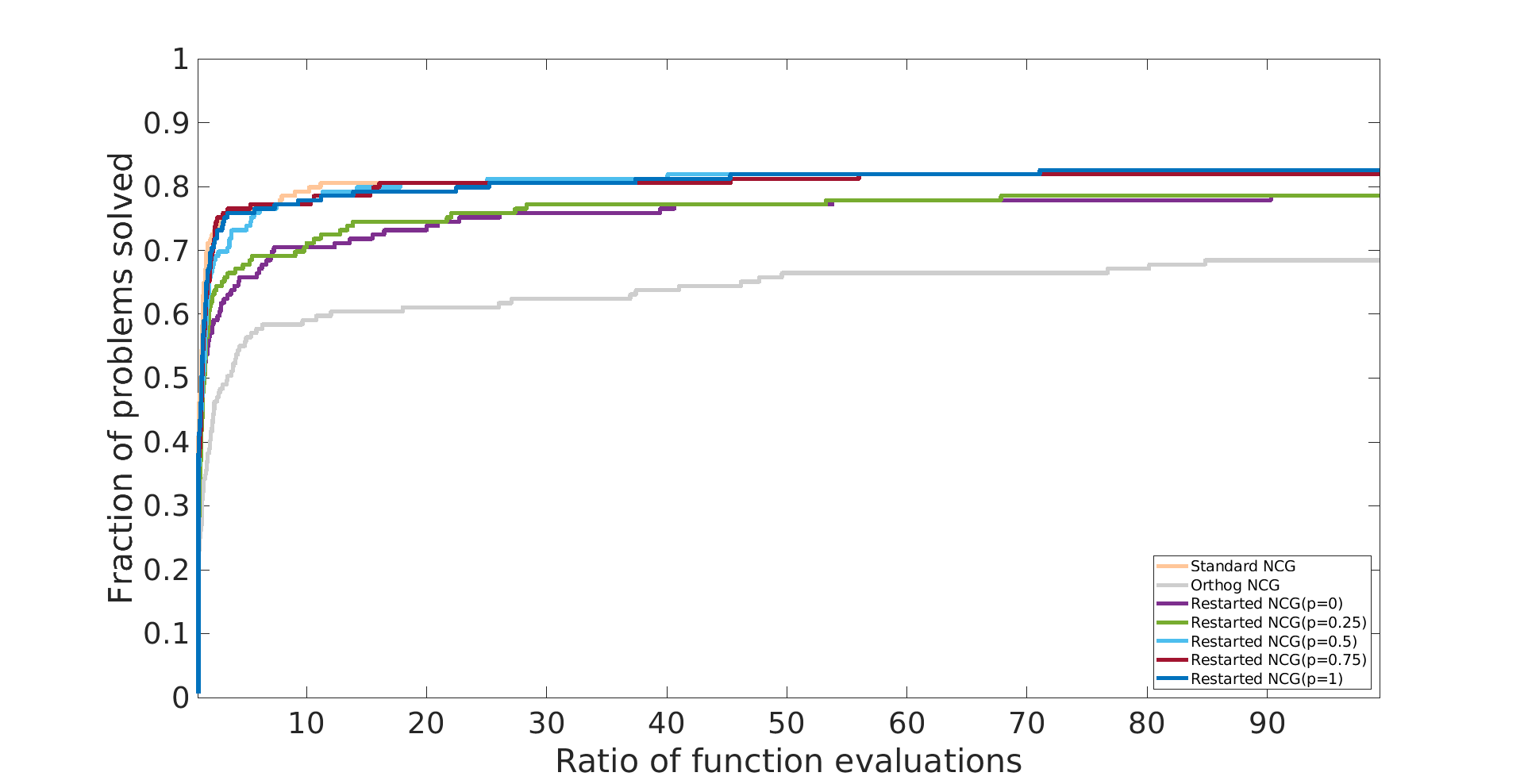}
  \caption{
  \revised{Performance profiles (function evaluations) 
  for nonlinear CG methods with the FR formula~\eqref{eq:ncgFR}
formula on a benchmark from CUTEst.}}
  \label{fig:cutest:FRfevals}
\end{figure}

\revised{We present our results under the form of performance 
profiles~\cite{EDDolan_JJMore_2002} for three choices of formula for 
$\beta_{k+1}$. We provide those profiles both using both the number of 
iterations (which would also correspond to the number of gradient 
evaluations) and the number of function evaluations as budget indicators.
}
\revised{When the FR formula~\eqref{eq:ncgFR} is used, we 
again observe that using $p \ge 0.5$ in Restarted NCG yields a profile 
quite close to that of Standard NCG. As illustrated by Figures~\ref{fig:cutest:FRits} 
and~\ref{fig:cutest:FRfevals}, the iteration and evaluation 
profiles are highly similar (yet not identical), suggesting that the 
number of backtracking line-search iterations is essentially constant.}

\begin{figure}
  \centering
  \includegraphics[width=\textwidth]{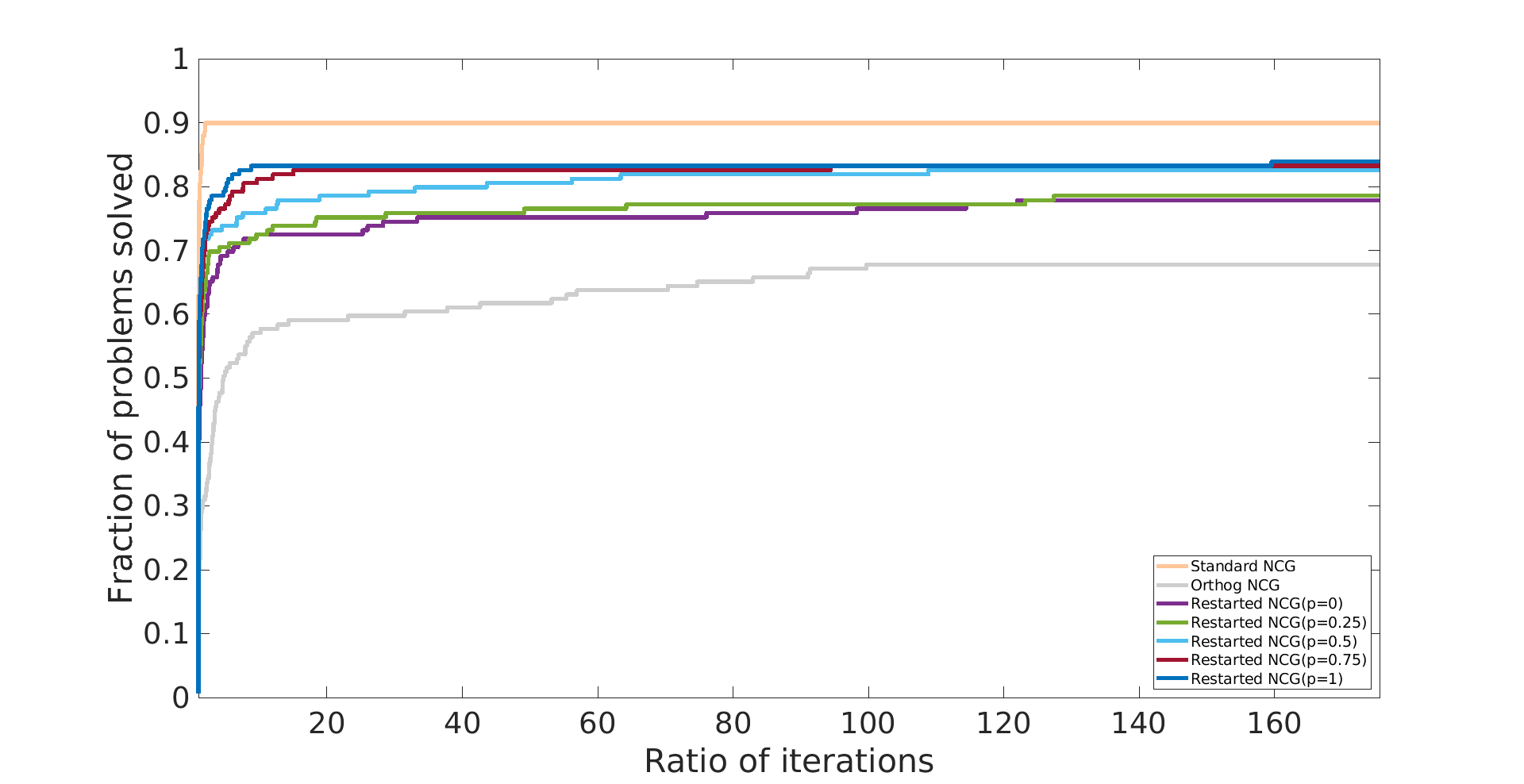}
  \caption{
  \revised{Performance profiles (iterations) 
  for nonlinear CG methods with the HZ formula~\eqref{eq:ncgHZ}
formula on a benchmark from CUTEst.}}
  \label{fig:cutest:HZits}
\end{figure}

\begin{figure}
  \centering
  \includegraphics[width=\textwidth]{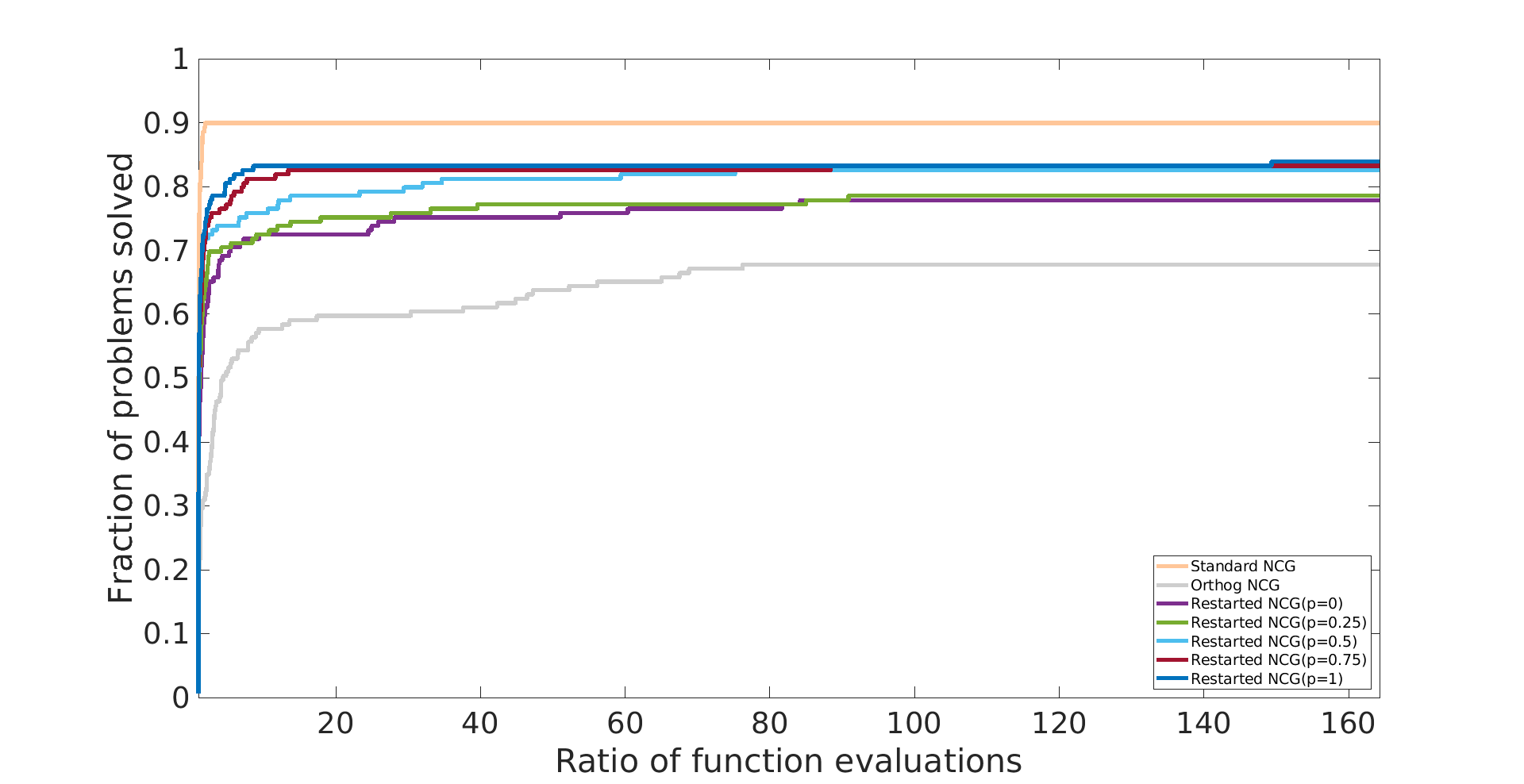}
  \caption{
  \revised{Performance profiles (function evaluations) 
  for nonlinear CG methods with the HZ formula~\eqref{eq:ncgHZ}
formula on a benchmark from CUTEst.}}
  \label{fig:cutest:HZfevals}
\end{figure}

\revised{Figures~\ref{fig:cutest:HZits} and~\ref{fig:cutest:HZfevals} 
show profiles obtained with the HZ formula~\eqref{eq:ncgHZ}. On these plots, Standard NCG is capable 
of solving a larger fraction of the problems than the Restarted 
NCG methods. However, we still observe that the best profile among 
the Restarted NCG methods are obtained for $p \ge 0.5$.}

\begin{figure}
  \centering
  \includegraphics[width=\textwidth]{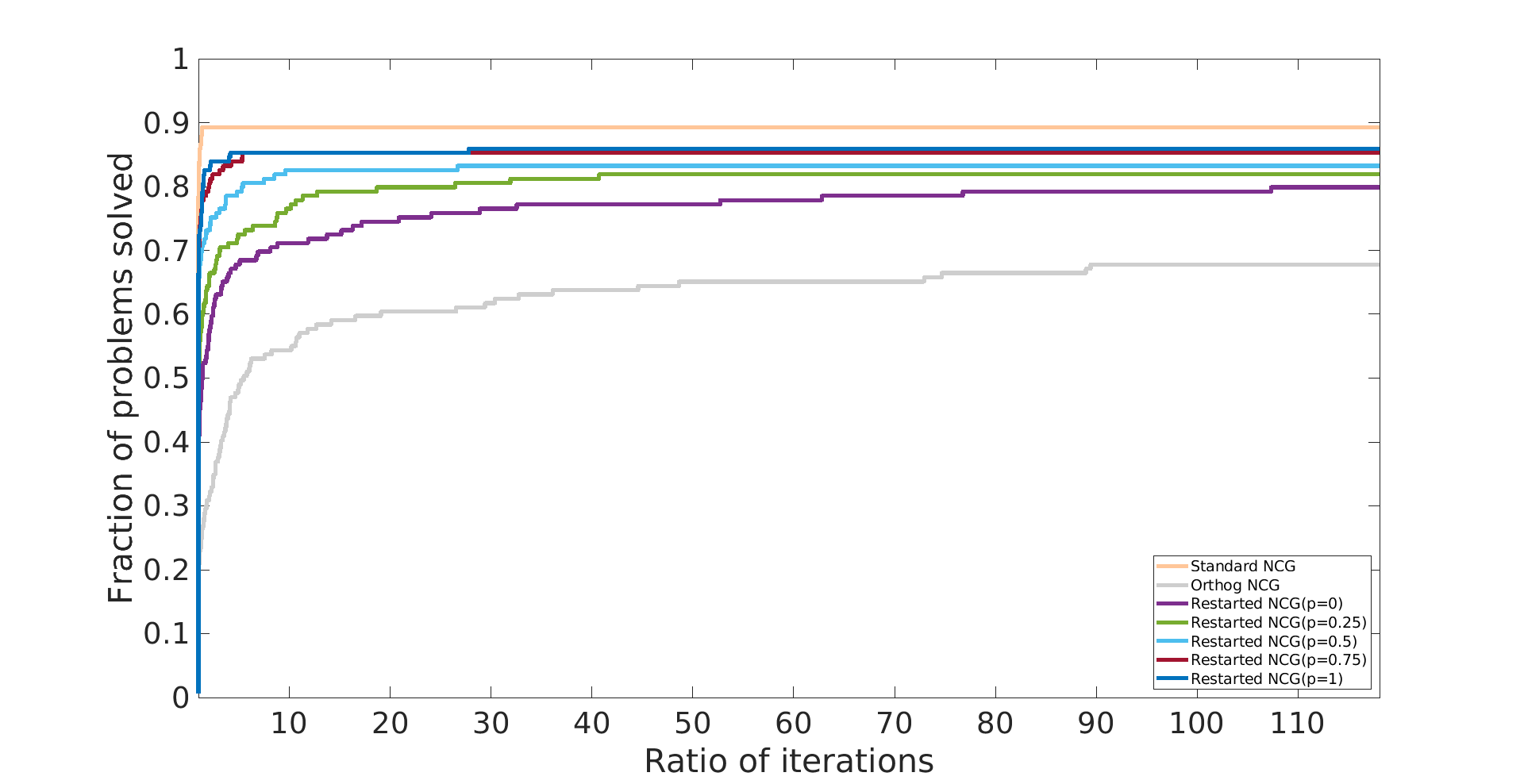}
  \caption{
  \revised{Performance profiles (iterations) 
  for nonlinear CG methods with the PRP+ formula~\eqref{eq:ncgPRP+} 
formula on a benchmark from CUTEst.}}
  \label{fig:cutest:PRPplusits}
\end{figure}

\begin{figure}
  \centering
  \includegraphics[width=\textwidth]{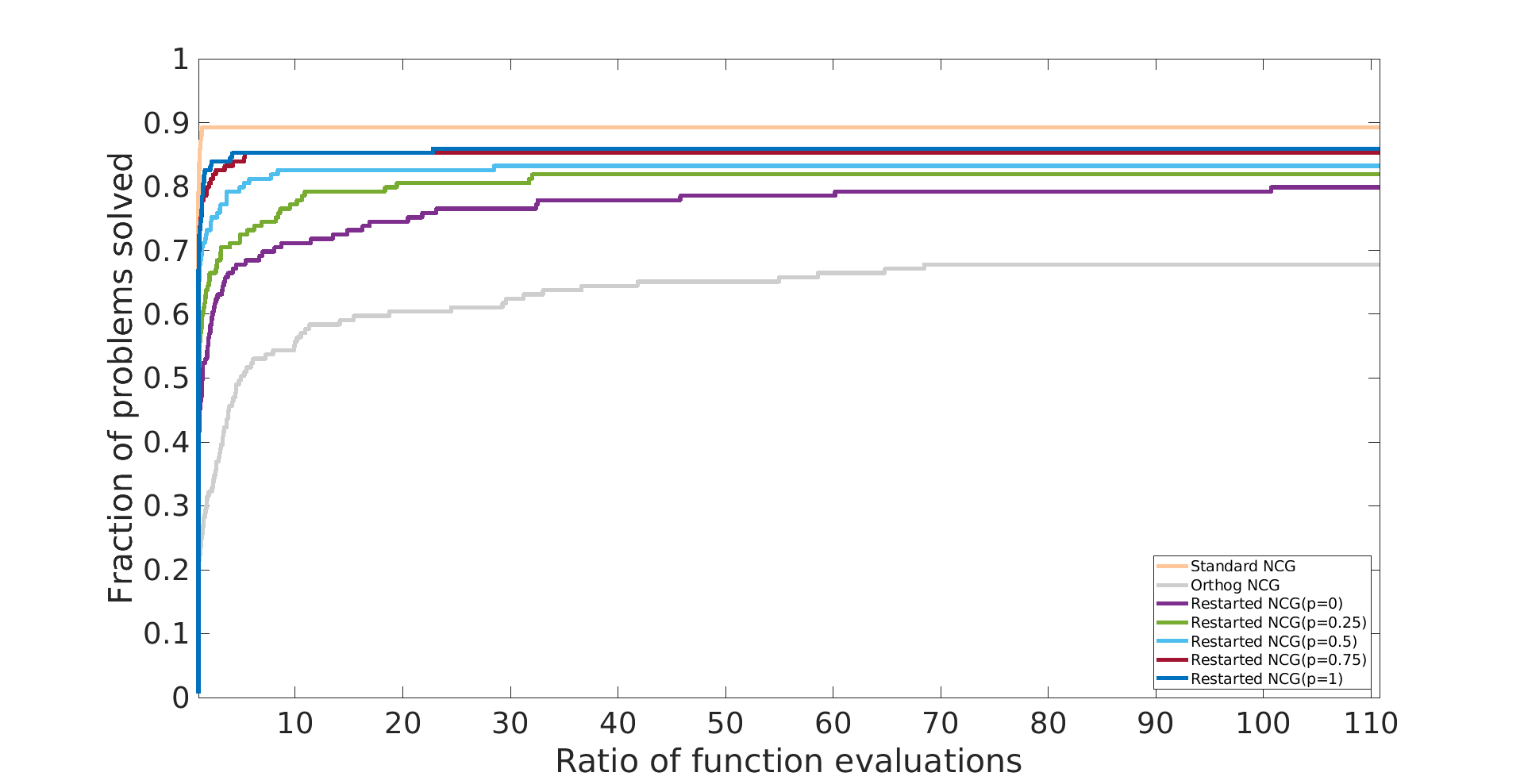}
  \caption{
  \revised{Performance profiles (function evaluations) 
  for nonlinear CG methods with the PRP+ formula~\eqref{eq:ncgPRP+}
formula on a benchmark from CUTEst.}}
  \label{fig:cutest:PRPplusfevals}
\end{figure}

\revised{Finally, we provide results obtained using the 
PRP+ formula~\eqref{eq:ncgPRP+} through 
Figures~\ref{fig:cutest:PRPplusits} and~\ref{fig:cutest:PRPplusfevals}. 
The results are of the same flavor than the previous ones, although 
the gap between Standard NCG and the best Restarted NCG methods is 
smaller.}

\revised{In a nutshell, these experiments suggest that the 
Restarted NCG variants may exhibit worse performance than a standard 
nonlinear conjugate gradient approach (that just restarts when it 
computes a non-descent search direction). The gap between a standard 
method and our restarted variant can be partially explained by the 
additional requirements put in place to ensure complexity guarantees. 
However, our experiments with the Fletcher-Reeves formula also show 
that the restarting condition can have quite a minor effect on the 
performance, which is encouraging for future investigation on these 
aspects. Finally, we note that the Orthog NCG variant did not perform 
well on this test set, regardless of the formula that was used. Our 
interpretation is that the loss of orthogonality here, as measured 
by the test $|g_k^\T g_{k+1}| \ge \sigma \|g_k\|^2$ is not 
necessarily detrimental to the algorithm performance. Though out of 
the scope of this work, we conjecture that varying the power within 
the right-hand side of the condition might improve the performance.}

%
%\begin{figure}
%\centering
%     \begin{subfigure}[b]{0.5\textwidth}
%         \centering
%         \includegraphics[width=1.4\textwidth]{./cutestFRits.png}
%         \caption{Iterations}
%         %\label{fig:y equals x}
%     \end{subfigure}
%     \hfill
%     \begin{subfigure}[b]{0.5\textwidth}
%         \centering
%         \includegraphics[width=1.4\textwidth]{./cutestFRfevals.png}
%         \caption{Objective calls}
%         %\label{fig:y equals x}
%     \end{subfigure}
%\caption{Performance profiles for nonlinear CG methods with $\beta_{k+1}^{FR}$ 
%formula.}
%\label{fig:cutest:FR}
%\end{figure}

%%%%%%%%%%%%%%%%%%%%%%%%%%%%%%%%%%%%%%%%%%%%%%%%%%%%%%%%%%%%%%%%%%%%%%%%%%%%%%%
\section{Conclusion}
\label{sec:conc}
%%%%%%%%%%%%%%%%%%%%%%%%%%%%%%%%%%%%%%%%%%%%%%%%%%%%%%%%%%%%%%%%%%%%%%%%%%%%%%%

In this paper, we presented a nonlinear conjugate gradient framework based on 
Armijo line search and a modified restart condition, which we endowed with
worst-case complexity guarantees. Although the results are of the same 
order than that for gradient descent, our complexity bound illustrates that 
better properties of nonlinear conjugate gradient may improve the overall 
number of iterations.
Our motivation for considering this particular framework was the remarkable 
performance of Armijo PRP+ nonlinear CG on a nonconvex regression problem. 
Our experiments \revised{on such instances} suggest that our new restart 
condition may be parameterized so as to match the original nonlinear CG 
method, \revised{thereby providing a theoretical justification of the 
performance of this method.}

Our analysis does not leverage the specific definition of the search 
directions, but rather checks their properties a posteriori. \revised{A 
natural continuation of the present paper would consist in enforcing 
similar properties by design, possibly by using stronger line-search 
conditions such as strong Wolfe. In addition, the design of a nonlinear 
conjugate gradient with strictly better complexity bounds than gradient 
descent remains an open question.} Recent advances in this
 area for 
convex optimization~\cite{SKarimi_SAVavasis_2021} might provide insights  
regarding the variants that are most amenable to a complexity analysis 
with guaranteed improvement over gradient descent.

%%%%%%%%%%%%%%%%%%%%%%%%%%%%%%%%%%%%%%%%%%%%%%%%%%%%%%%%%%%%%%%%%%%%%%%%%%%%%%%
%\nocite{*}
\bibliographystyle{plain}
\bibliography{biblioncg}
%%%%%%%%%%%%%%%%%%%%%%%%%%%%%%%%%%%%%%%%%%%%%%%%%%%%%%%%%%%%%%%%%%%%%%%%%%%%%%%

\end{document}